\documentclass[11pt]{amsart}
\usepackage[cp1251]{inputenc}

\newtheorem{theorem}{Theorem}[section]
\newtheorem{lemma}[theorem]{Lemma}
\newtheorem{proposition}[theorem]{Proposition}
\newtheorem{corollary}[theorem]{Corollary}

\theoremstyle{definition}
\newtheorem{definition}[theorem]{Definition}

\newcommand{\Aut}{\operatorname{Aut}}
\newcommand{\TAut}{\operatorname{TAut}}
\newcommand{\NAut}{\operatorname{NAut}}

\def\Ind{\operatorname{Ind}}
\def\SL{\operatorname{SL}}
\def\End{\operatorname{End}}
\def\Sympl{\operatorname{Sympl}}
\def\Id{\operatorname{Id}}
\def\Ch{\operatorname{Char}}
\def\Fr{\operatorname{Fr}}

\theoremstyle{remark}
\usepackage{amscd,amssymb}
\begin{document}
\renewcommand{\to}{\mapsto}
\vskip-5cm

\baselineskip 19 pt

\title[Tame approximation] {On the tame authomorphism approximation, augmentation Topology of
Automorphism Groups and $Ind$-schemes, and authomorphisms of tame automorphism groups}
\voffset-2cm

\author[A. Kanel-Belov, J.-T. Yu, and A. Elishev]
{Alexei Kanel-Belov, Jie-Tai Yu, and Andrey Elishev}
\address{College of Mathematics and Statistics, Shenzhen University, Shenzhen, 518061, China} \email{jietaiyu@szu.edu.cn}, \email{ kanelster@gmail.com}
\address{Laboratory of Advanced Combinatorics and Network Applications, Moscow Institute of Physics and Technology, Dolgoprudny, Moscow Region, 141700, Russia}
\email{elishev@phystech.edu}
\thanks{The research was supported by Russian Science Foundation, Grant No 17-11-01377}.




\begin{abstract}

We study topological properties of $\Ind$-groups
$\Aut(K[x_1,\dots,x_n])$ and \\
$\Aut(K\langle x_1,\dots,x_n\rangle)$ of automorphisms of polynomial and free
associative algebras
via $\Ind$-schemes, toric varieties,
approximations, and singularities.

We obtain a number of properties of $\Aut(\Aut(A))$, where $A$ is the
polynomial or free associative algebra over the base field $K$. We prove
that all $\Ind$-scheme automorphisms of $\Aut(K[x_1,\dots,x_n])$ are
inner for $n\ge 3$, and all $\Ind$-scheme automorphisms of\\
$\Aut(K\langle x_1,\dots, x_n\rangle)$ are semi-inner.

As an application, we prove that $\Aut(K[x_1,\dots,x_n])$
cannot be embedded into \\
$\Aut(K\langle x_1,\dots,x_n\rangle)$
by the natural abelianization.
In other words, the {\it Automorphism Group Lifting Problem}
has a negative solution.

We explore close connection between the above results and the
Jacobian conjecture, as well as the Kanel-Belov -- Kontsevich conjecture, and
formulate the Jacobian conjecture for fields of any
characteristic.

We make use of results of Bodnarchuk and Rips, and we also consider
automorphisms of tame groups preserving the origin and obtain a modification of said results in the tame setting.

%
%
%
%
%
%
%
%

\end{abstract}

\maketitle

\qquad \qquad Dedicated to prof. B.I.Plotkin who initiated the subject

\tableofcontents

\section{Introduction and main results}

\subsection{Automorphisms of $K[x_1,\dots,x_n]$ and $K\langle x_1,\dots,x_n\rangle$}

Let $K$ be a field. The main objects of this study are the $K$-algebra automorphism groups $\Aut K[x_1,\dots,x_n]$ and \\
$\Aut K\langle x_1,\dots,x_n\rangle$ of the (commutative) polynomial algebra and the free associative algebra with $n$ generators, respectively. The former is equivalent to the group of all polynomial one-to-one mappings of the affine space $\mathbb{A}^n_{K}$. Both groups admit a representation as a colimit of algebraic sets of automorphisms filtered by total degree (with morphisms in the direct system given by closed embeddings) which turns them into topological spaces with Zariski topology compatible with the group structure. The two groups carry a power series topology as well, since every automorphism $\varphi$ may be identified with the $n$-tuple $(\varphi(x_1),\ldots,\varphi(x_n))$ of the images of generators. This topology plays an especially important role in the applications, and it turns out -- as reflected in the main results of this study -- that approximation properties arising from this topology agree well with properties of combinatorial nature.

$\Ind$-groups of polynomial automorphisms play a central part in the study of the Jacobian
conjecture of O. Keller as well as a number of problems of similar nature. One outstanding example is provided by a recent conjecture of Kanel-Belov and Kontsevich (B-KKC), \cite{BelovKontsevich1,BelovKontsevich}, which asks whether
the group
$$
\Sympl({\mathbb C}^{2n})\subset \Aut(\mathbb{C}[x_1,\dots,x_{2n}])
$$
of complex polynomial automorphisms preserving the standard Poisson bracket $$\lbrace x_i,\;x_j\rbrace = \delta_{i,n+j}-\delta_{i+n,j}$$ is isomorphic\footnote{In fact, the conjecture seeks to establish an isomorphism
 $\Sympl(K^{2n})\simeq \Aut(W_n(K))$ for any field $K$ of characteristic zero in a functorial manner.} to the group of automorphisms of the $n$-th Weyl algebra $W_n$
\begin{gather*}
W_n(\mathbb{C})=\mathbb{C}\langle x_1,\ldots,x_n,y_1,\ldots,y_n\rangle/I,\\
I= \left( x_ix_j-x_jx_i,\; y_iy_j-y_jy_i,\; y_ix_j-x_jy_i-\delta_{ij}\right).
\end{gather*}

The physical meaning of Kanel-Belov and Kontsevich conjecture is the invariance of the polynomial symplectomorphism group of the phase space under the procedure of deformation quantization.

The B-KKC was conceived  during a successful search for a proof of stable equivalence
of the Jacobian conjecture and a well-known conjecture of Dixmier stating that
$\Aut(W_n)=\End(W_n)$ over any field of characteristic zero. In the papers \cite{BelovKontsevich1,BelovKontsevich}
a particular family of homomorphisms (in effect, monomorphisms) $\Aut(W_n(\mathbb{C}))\rightarrow\Sympl({\mathbb
C}^{2n})$ was constructed, and a natural question whether those homomorphisms were in fact isomorphisms was raised.
The aforementioned morphisms, independently studied by Tsuchimoto to the same end, were in actuality defined as restrictions of morphisms of the saturated model of Weyl algebra over an algebraically closed field of {\it positive} characteristic - an object which contains $W_n(\mathbb{C})$ as a proper subalgebra. One of the defined morphisms turned out to have a particularly simple form over the subgroup of the so-called tame automorphisms, and it was natural to assume that morphism was the desired B-KK isomorphism (at least for the case of algebraically closed base field). Central to the construction is the notion of infinitely large prime number (in the sense of hyperintegers), which arises as the sequence $(p_m)_{m\in\mathbb{N}}$ of positive characteristics of finite fields comprising the saturated model.  This leads to the natural problem
(\cite{BelovKontsevich}):

\medskip
\noindent {\bf Problem.} Prove that the B-KK morphism is independent of the choice of the infinite prime $(p_m)_{m\in\mathbb{N}}$.
\medskip

A general formulation of this question in the paper
\cite{BelovKontsevich} goes as follows:

For a commutative ring $R$ define
\begin{equation*}
R_{\infty}=\lim_{\rightarrow}\left( \prod_{p} R'\otimes \mathbb{Z}/p\mathbb{Z}\;/\;\bigoplus_{p} R'\otimes \mathbb{Z}/p\mathbb{Z}\right),
\end{equation*}
where the direct limit is taken over the filtered system of all finitely generated subrings $R'\subset R$ and the product and the sum are taken over all primes $p$. This larger ring possesses a unique "nonstandard Frobenius" endomorphism $\Fr:R_{\infty}\rightarrow R_{\infty}$ given by
\begin{equation*}
(a_p)_{\text{primes\;}p}\mapsto (a_p^p)_{\text{primes\;}p}.
\end{equation*}

The Kanel-Belov and Kontsevich construction returns a morphism
\begin{equation*}
\psi_R: \Aut (W_n(R))\rightarrow \Sympl R_{\infty}^{2n}
\end{equation*}

such that there exists a unique homomorphism
$$\phi_R: \Aut(W_n)(R)\rightarrow \Aut(P_n)(R_\infty)$$
obeying $\psi_R=\Fr_*\circ \phi_R$. Here
$\Fr_*:\Aut(P_n)(R_\infty)\rightarrow \Aut(P_n)(R_\infty)$
 is the $\Ind$-group homomorphism induced by the Frobenius endomorphism of the coefficient ring, and $P_n$ is the commutative Poisson algebra, i.e. the polynomial algebra in $2n$ variables equipped with additional Poisson structure (so that $\Aut( P_n(R))$ is just $\Sympl( R^{2n})$ - the group of Poisson structure-preserving automorphisms).

\medskip
\noindent {\bf Question.}  {\it In the above formulation, does the image of $\phi_R$
belong to
$$\Aut(P_n)(i(R)\otimes {\mathbb Q})\,\,,$$
where $i:R\rightarrow R_\infty$ is the tautological inclusion? In other
words, does there exist a unique homomorphism
$$\phi_R^{can}:\Aut(P_n)(R)\rightarrow \Aut(P_n)(R\otimes {\mathbb Q})$$
 such that $\psi_R=\Fr_*\circ i_*\circ \phi_R^{can}$.}

\medskip

Comparing the two morphisms $\phi$ and $\varphi$ defined using two
different free ultrafilters, we obtain a "loop" element $\phi\varphi^{-1}$ of
 $\Aut_{\Ind}(\Aut(W_n))$,
(i.e. an automorphism which preserves the structure of infinite dimensional
algebraic group). Describing this group would provide a solution to this
question.

Some progress toward resolution of the B-KKC independence problem has been made recently in \cite{KBE1, KBE2}, although the general unconditional case is still open.

\bigskip

In the spirit of the above we propose the following

\medskip
\noindent {\bf Conjecture.}
{\it All automorphisms of the $\Ind$-group $\Sympl({\mathbb
C}^{2n})$ are inner.}

A similar conjecture may be put forward for $\Aut(W_n(\mathbb{C}))$.

We are focused on the  investigation of the group
$\Aut(\Aut(K[x_1,\dots,x_n]))$ and the corresponding noncommutative (free associative
algebra) case. This way of thinking has its roots in the realm of universal algebra and universal algebraic geometry and was conceived in the pioneering work of Boris Plotkin. A more detailed discussion can be found in \cite{BelovLipBer}.



\

\noindent {\bf Wild automorphisms and the lifting problem.} In 2004, the celebrated Nagata
conjecture over a field $K$ of characteristic zero was proved by
Shestakov and Umirbaev \cite{SU1, SU2} and a stronger version of
the conjecture was proved by Umirbaev and Yu \cite{UY}. Let $K$ be a field of characteristic zero. Every wild
$K[z]$-automorphism (wild $K[z]$-coordinate) of $K[z][x,y]$ is
wild viewed as a $K$-automorphism ($K$-coordinate) of $K[x,y,z]$.
In particular, the Nagata automorphism $(x-2y(y^2+xz)-(y^2+xz)^2z,
y+(y^2+xz)z, z)$ (Nagata coordinates $x-2y(y^2+xz)-(y^2+xz)^2z$
and $y+(y^2+xz)z$) are wild. In \cite{UY}, a related question
was raised:

\

\noindent {\bf The lifting problem.} \emph{Can an arbitrary wild
automorphism (wild coordinate) of the polynomial algebra
$K[x,y,z]$ over a field $K$ be lifted to an automorphism
(coordinate) of the free associative algebra $K\langle x,y,z\rangle$?}

\medskip

In the paper \cite{BelovYuLifting}, based on the degree estimate
\cite{MLY, YuYungChang}, it was proved that any wild
$z$-automorphism including the Nagata automorphism cannot be
lifted as a $z$-automorphism (moreover, in \cite{BKY}
it is proved that every $z$-automorphism of $K\langle x, y,z\rangle$ is
stably tame and becomes tame after adding at most one variable).
It means that  if every automorphism can be lifted, then it provides
an obstruction $z'$ to $z$-lifting and the question to estimate  such an
obstruction is naturally raised.

In view of the above, we may ask the following:

\

 \noindent {\bf The automorphism group lifting problem.}
\emph{Is $\Aut(K[x_1,\dots,x_n])$ isomorphic to a subgroup of
$\Aut(K\langle x_1,\dots,x_n\rangle)$ under the natural
abelianization?}

\medskip

The following examples show this problem is interesting and
non-trivial.

\medskip
\noindent{\bf Example 1.} There is a surjective
homomorphism (taking the absolute value) from $\mathbb C^*$ onto $\mathbb R^+$.
But $\mathbb R^+$ is isomorphic to the subgroup $\mathbb R^+$ of $\mathbb C^*$
under the homomorphism.

\medskip

\noindent{\bf Example 2.} There is a surjective homomorphism
(taking the determinant) from $\text{GL}_n(\mathbb R)$ onto
$\mathbb R^*$. But obviously $\mathbb R^*$ is isomorphic to the
subgroup $\mathbb R^*I_n$ of $\text{GL}_n(\mathbb R)$.

\medskip

In this paper we prove that the automorphism group lifting problem has a negative answer.

The lifting problem  and the automorphism group lifting problem are closely related to the Kanel-Belov and Kontsevich
Conjecture (see Section \ref{KBConj}).

Consider a symplectomorphism $\varphi: x_i\to P_i,\; y_i\to Q_i$. It
can be lifted to some automorphism $\widehat{\varphi}$ of the
quantized algebra $W_\hbar [[\hbar]]$: $$\widehat{\varphi}: x_i\to
P_i+P_i^1\hbar+\cdots+P_i^m\hbar^m;\ y_i\to
Q_i+Q_i^1\hbar+\cdots+Q_i^m\hbar^m.$$ The point is to choose a lift
$\widehat{\varphi}$ in such a way that the degree of all $P_i^m,
Q_i^m$ would be bounded. If that is true, then the B-KKC
follows.
\bigskip

\subsection{Main results}
The main results of this paper are as follows.

 \begin{theorem}        \label{ThAutAut}
Any $\Ind$-scheme automorphism $\varphi$ of
$\NAut(K[x_1,\dots,x_n])$ for $n\ge 3$ is inner, i.e. is a
conjugation via some automorphism.
\end{theorem}

 \begin{theorem}        \label{ThAutAutFreeAssoc}
Any $\Ind$-scheme automorphism $\varphi$ of
$\NAut(K\langle
x_1,\dots,x_n\rangle)$ for $n\ge 3$ is semi-inner (see definition
\ref{DfSemi}).
\end{theorem}

$NAut$ denotes the group of {\em nice} automorphisms, i.e. automorphisms
which can be approximated by tame ones (definition
\ref{DfNiceAprox}). In characteristic zero case every automorphism
is nice.

For the group of  automorphisms of a semigroup a number of similar results
on set-theoretical level was obtained previously by Kanel-Belov,
Lipyanski and Berzinsh \cite{BelovLiapiansk2, BelovLipBer}.
All these questions (including $\Aut(\Aut)$ investigation) take root in the realm of Universal Algebraic Geometry and were proposed
by Boris Plotkin. Equivalence of two algebras having the same generalized
identities and isomorphism of first order means semi-inner
properties of automorphisms (see
\cite{BelovLiapiansk2, BelovLipBer} for details).

\

\noindent {\bf Automorphisms of tame automorphism groups.}
Regarding the tame automorphism group, something can be done on
the group-
theoretic level. In the paper of H. Kraft and
I. Stampfli \cite{KraftStampfli} the automorphism group of the tame
automorphism group of the polynomial algebra was thoroughly studied.
In that paper, conjugation of elementary automorphisms via
translations played an important role. The results of our study are different. We describe the group $\Aut(\TAut_0)$ of the group $\TAut_0$ of tame automorphisms preserving the origin (i.e. taking
the augmentation ideal onto an ideal which is a subset of the
augmentation ideal). This is technically more difficult, and
will be universally and systematically done for both commutative
(polynomial algebra) case and noncommutative (free associative
algebra) case. We observe a few problems in the shift conjugation
approach for the noncommutative (free associative algebra) case,
as it was for commutative case in \cite{KraftStampfli}. Any
evaluation on a ground field element can return zero,
for example in Lie polynomial $[[x,y],z]$. Note that the calculations
of $\Aut(\TAut_0)$ (resp. $\Aut_{\Ind}(\TAut_0)$,
$\Aut_{\Ind}(\Aut_0)$) imply also the same results for
$\Aut(\TAut)$ (resp. $\Aut_{\Ind}(\TAut)$, $\Aut_{\Ind}(\Aut)$)
according to the approach of this article via stabilization by the
torus action.

\begin{theorem}        \label{ThAutTAut}
Any  automorphism $\varphi$ of $\TAut_0(K[x_1,\dots,x_n])$ (in the
group-theoretic sense) for $n\ge 3$ is inner, i.e. is a
conjugation via some automorphism.
\end{theorem}

 \begin{theorem}        \label{ThAutTAut0}
The group $\TAut_0(K[x_1,\dots,x_n])$ is generated by the
automorphism $$x_1\to x_1+ x_2x_3,\; x_i\to x_i, \;\;i\ne 1$$ and linear substitutions
if $\Ch(K)\ne 2$ and $n>3$.
\end{theorem}

Let $G_N\subset \TAut(K[x_1,\dots,x_n])$, $E_N\subset \TAut(K\langle
x_1,\dots,x_n\rangle)$ be tame automorphism subgroups preserving
the $N$-th power of the augmentation ideal.

 \begin{theorem}        \label{ThAutTAutOr}
Any  automorphism $\varphi$ of $G_N$ (in the group-theoretic sense)
for $N\ge 3$ is inner, i.e. is given by a conjugation via some
automorphism.
\end{theorem}

\begin{definition}   \label{DfSemi}
An {\em anti-automorphism} $\Psi$ of a $K$-algebra $B$ is a vector space automorphism such that $\Psi(ab)=\Psi(b)\Psi(a)$.
For instance, transposition of matrices is an anti-automorphism.
An anti-automorphism of the free associative algebra $A$ is a {\em mirror anti-automorphism}
if it sends $x_ix_j$ to $x_jx_i$ for some fixed $i$ and $j$. If a mirror anti-automorphism $\theta$ acts identical on all generators $x_i$, then for any monomial $x_{i_1}\cdots x_{i_k}$ we have
$$
\theta(x_{i_1}\cdots x_{i_k})=x_{i_k}\cdots x_{i_1}.
$$
Such an anti-automorphism will be generally referred to as \emph{the} mirror anti-automorphism.

An automorphism of $\Aut(A)$ is {\em semi-inner} if it can be
expressed as a composition of an inner automorphism and a conjugation by a mirror anti-automorphism.
\end{definition}

\begin{theorem}        \label{ThAutTAutFreass}
a) Any  automorphism $\varphi$ of $\TAut_0(K\langle
x_1,\dots,x_n\rangle)$ and also\\
 $\TAut(K\langle
x_1,\dots,x_n\rangle)$ (in the group-theoretic sense) for $n\ge 4$
is semi-inner, i.e. is a conjugation via some automorphism and/or
mirror anti-automorphism.

b) The same is true for $E_n$, $n\ge 4$.
\end{theorem}

The case of $\TAut(K\langle x,y,z\rangle)$ is substantially more difficult.
We can treat it only on $\Ind$-scheme level, but even then it is
the most technical part of the paper (see section
\ref{SbSc3VrbFreeAss}). For the two-variable case a similar proposition is probably false.

\begin{theorem}            \label{ThTAss3Ind}
a) Let $\Ch(K)\ne 2$. Then $\Aut_{\Ind}(\TAut(K\langle
x,y,z\rangle))$ (resp. \\
$\Aut_{\Ind}(\TAut_0(K\langle x,y,z\rangle))$) is generated by conjugation by an automorphism or
a mirror anti-automorphism.

b) The same is true for $\Aut_{\Ind}(E_3)$.
\end{theorem}

By $\TAut$ we denote the tame automorphism group, $\Aut_{\Ind}$ is the
group of $\Ind$-scheme automorphisms (see section
\ref{ScIndShme}).

Approximation allows us to formulate  the celebrated Jacobian conjecture for
any characteristic.

\

{\bf Lifting of the automorphism groups.} In this article we prove
that the automorphism group of
 polynomial algebra over an arbitrary field $K$ cannot be embedded
 into the automorphism group of free associative algebra induced by the  natural abelianization.

\begin{theorem}     \label{ThGroupLifting}
Let $K$ be an arbitrary field, $G=\Aut_0(K[x_1,\dots,x_n])$ and
$n>2$. Then $G$ cannot be isomorphic to any subgroup $H$ of
$\Aut(K\langle x_1,\dots,x_n\rangle)$ induced by the natural
abelianization. The same is true for $\NAut(K[x_1,\dots,x_n])$.
\end{theorem}


\section{Varieties of automorphisms}

\subsection{Elementary and tame automorphisms}

Let $P$ be a polynomial that is independent of $x_i$ with $i$ fixed. An automorphism
$$x_i\to x_i+P,\; x_j\to x_j\;\; \mbox{for}\ i\ne j$$ is called {\em
elementary}. The group generated by linear automorphisms and
elementary ones for all possible $P$ is called the {\em tame
automorphism group (or subgroup) $\TAut$} and elements of $\TAut$ are {\em tame
automorphisms}.

\subsection{$\Ind$-schemes and $\Ind$-groups} \label{ScIndShme}

\begin{definition}
An {\it $\Ind$-variety} $M$ is the direct limit of algebraic varieties
$M=\varinjlim \lbrace M_1\subseteq M_2\cdots\rbrace$.
 An {\it $\Ind$-scheme} is an
$\Ind$-variety which is a group such that the group inversion is
a morphism $M_i\rightarrow M_{j(i)}$ of algebraic varieties, and the group multiplication induces
a morphism from $M_i\times M_j$ to $M_{k(i,j)}$. A map $\varphi$ is a {\it
morphism} of an $\Ind$-variety $M$ to an $\Ind$-variety $N$, if
$\varphi(M_i)\subseteq N_{j(i)}$ and the restriction $\varphi$ to
$M_i$ is a morphism for all $i$. Monomorphisms, epimorphisms and
isomorphisms are defined similarly in a natural way.
\end{definition}

{\bf Example.} $M$ is the group of automorphisms of the affine space, and
$M_j$ are the sets of all automorphisms in $M$ with degree $\le j$.
\medskip

There is an interesting

\medskip
{\bf Problem.} \emph{Investigate growth functions of
$\Ind$-varieties. For example, the dimension of varieties of
polynomial automorphisms of degree $\le n$.}

\medskip

Note that coincidence of growth functions of $\Aut(W_n(\mathbb{C}))$ and
$\Sympl({\mathbb C}^{2n})$ would imply the Kanel-Belov -- Kontsevich conjecture
\cite{BelovKontsevich}.

\begin{definition}             \label{DfAugm}
The ideal $I$ generated by variables $x_i$ is called the {\em
augmentation ideal}.  For a fixed positive integer $N>1$, the {\em augmentation subgroup $H_N$} is the
group of all automorphisms $\varphi$ such that $\varphi(x_i)\equiv
x_i \mod I^N$. The larger group $\hat{H}_N\supset H_N$ is the group of
automorphisms whose linear part is scalar, and $\varphi(x_i)\equiv
\lambda x_i \mod I^N$ ($\lambda$ does not depend on $i$). We often say an arbitrary element of the group $\hat{H}_N$ is an automorphism that is homothety modulo (the $N$-th power of) the augmentation ideal.
\end{definition}

\section{The Jacobian conjecture in any characteristic,
Kanel-Belov -- Kontsevich conjecture, and approximation}

\subsection{Approximation problems and Kanel-Belov -- Kontsevich Conjecture}    \label{KBConj}

Let us give formulation of the Kanel-Belov -- Kontsevich Conjecture:

\medskip
{\bf $B-KKC_n$}: $\Aut(W_n)\simeq\Sympl({\mathbb C}^{2n})$.
\medskip

A similar conjecture can be stated for endomorphisms

\medskip
{\bf $B-KKC_n$}: $\End(W_n)\simeq\Sympl\End({\mathbb C}^{2n})$.
\medskip

If the Jacobian conjecture $JC_{2n}$ is true, then the respective conjunctions over all $n$ of the two
conjectures are equivalent. 

It is natural to  approximate  automorphisms by tame ones. There
exists such an approximation up to terms of any order for polynomial automorphisms as well as Weyl algebra automorphisms, symplectomorphisms etc. However, the naive approach
fails.


It is known  that $\Aut(W_1)\equiv \Aut_1(K[x,y])$ where $\Aut_1$
stands for the subgroup of automorphisms of Jacobian determinant one. However, considerations
from \cite{Shafarevich} show that Lie algebra of the first group
is the algebra of derivations of $W_1$ and thus possesses no identities apart from the ones
of the free Lie algebra, another coincidence of the vector fields
which diverge to zero, and has polynomial identities. These cannot be
isomorphic \cite{BelovKontsevich1,BelovKontsevich}. In other
words, this group has two coordinate system  non-smooth  with
respect to one another (but integral with respect to one another). One
system is built from the coefficients of differential operators in a fixed basis of generators, while its counterpart is
provided by the coefficients of polynomials, which are images of the basis
$\tilde{x}_i, \tilde{y}_i$. 

In the paper \cite{Shafarevich} functionals on ${\mathfrak
m}/{\mathfrak m}^2$ were considered in order to define the Lie
algebra structure. In the spirit of that we have the following

\medskip
{\bf Conjecture.} The natural limit of ${\mathfrak m}/{\mathfrak m}^2$
is zero.
\medskip

It means that the definition of  the Lie algebra admits some sort of
functoriality problem and it depends on the presentation of
(reducible) $\Ind$-scheme.

In his remarkable paper,  Yu. Bodnarchuk \cite{Bodnarchuk}
established Theorem \ref{ThAutAut} by using Shafarevich's
results for the tame automorphism subgroup and for the case when the $\Ind$-scheme automorphism
is regular in the sense that it sends coordinate functions to coordinate functions. In this case the
tame approximation works (as well as for the symplectic case), and the corresponding method is similar to ours. We present
it here in order to make the text more self-contained, as well as
for the purpose of tackling the noncommutative (that is, the free associative algebra) case.
Note that in general, for regular functions, if the Shafarevich-style approximation were valid, then the Kanel-Belov -- Kontsevich
conjecture
would follow directly, which is absurd.


In the sequel, we do not assume regularity in the sense of \cite{Bodnarchuk} but
only assume that the restriction of a morphism on any subvariety is a morphism again. Note
that morphisms of $\Ind$-schemes $\Aut(W_n)\rightarrow\Sympl({\mathbb
C}^{2n})$ have this property, but are not regular in the sense of
Bodnarchuk \cite{Bodnarchuk}.

We use the idea of
singularity  which allows us to prove the augmentation subgroup
structure preservation, so that the approximation works in this case.

Consider the isomorphism $\Aut(W_1)\cong \Aut_1(K[x,y])$. It has
a strange property. Let us add a small parameter $t$. Then an
element arbitrary close to zero with respect to $t^k$ does not go
to zero  arbitrarily, so it is impossible to make tame limit!
There is a sequence of convergent product of elementary
automorphisms, which is not convergent under this isomorphism.
Exactly the same situation happens  for $W_n$. These effects cause
problems in perturbative quantum field theory.

\medskip

\subsection{The Jacobian conjecture in any characteristic}
Recall that the Jacobian conjecture in characteristic zero states that any polynomial endomorphism
\begin{equation*}
\varphi:K^n\rightarrow K^n
\end{equation*}
with constant Jacobian is globally invertible.

A naive attempt to directly transfer this formulation to positive characteristic fails because of the counterexample $x\mapsto
x-x^p$ ($p=\Ch K$), whose Jacobian is everywhere $1$ but which is evidently not invertible. Approximation provides a way to
formulate a suitable generalization of the Jacobian conjecture to any
characteristic and put it in a framework of other questions.


\begin{definition}        \label{DfNiceAprox}
An endomorphism $\varphi\in\End(K[x_1,\dots,x_n])$ is {\em good} if\\
for any $m$ there exist $\psi_m\in\End(K[x_1,\dots,x_n])$ and\\
$\phi_m\in\Aut(K[x_1,\dots,x_n])$ such that
\begin{itemize}
    \item $\varphi=\psi_m\phi_m$
    \item $\psi_m(x_i)\equiv x_i\mod (x_1,\dots,x_n)^m$.
\end{itemize}

An automorphism $\varphi\in\Aut(K[x_1,\dots,x_n])$ is {\em nice} if
for any $m$ there exist $\psi_m\in\Aut(K[x_1,\dots,x_n])$ and
$\phi_m\in\TAut(K[x_1,\dots,x_n])$ such that
\begin{itemize}
    \item $\varphi=\psi_m\phi_m$
    \item $\psi_m(x_i)\equiv x_i\mod (x_1,\dots,x_n)^m$, i.e. $\psi_m\in H_m$.
\end{itemize}
\end{definition}

Anick \cite{Anick} has shown that if $\Ch(K)=0$, any automorphism
is nice. However, this is unclear in positive characteristic.

\

{\bf Question.}\ {\it Is any automorphism over arbitrary field
nice?}
\medskip

Ever good automorphism has Jacobian $1$, and all such
automorphisms are good - and even nice - when $\Ch(K)=0$. This observation allows for the following question to be considered a generalization of the Jacobian conjecture to positive characteristic.

\medskip
{\bf The Jacobian conjecture in any characteristic:}\ {\it Is any good
endomorphism over arbitrary field an automorphism?}
\medskip

Similar notions can be formulated for the free associative algebra. That justifies the following

\medskip
{\bf Question.}\ {\it Is any automorphism of free associative
algebra over arbitrary field nice?}
\medskip

\

{\bf Question (version of free associative positive characteristic case of JC).}\ {\it Is any good endomorphism of the free associative
algebra over arbitrary field an automorphism?}
\medskip

\subsection{Approximation for the  automorphism group of
 affine spaces}
\label{SbScAprox}

\medskip

Approximation is the most important tool utilized in this paper. In order to perform it,
we have to prove that
$\varphi\in
\Aut_{\Ind}(\Aut_0(K[x_1,\dots,x_n])$ preserves the structure of
the augmentation subgroup.

The proof method utilized in theorems below works for commutative associative and free associative case. It is a problem of considerable interest to develop similar statements for automorphisms of other associative algebras, such as the commutative Poisson algebra (for which the $\Aut$ functor returns the group of polynomial symplectomorphisms); however, the situation there is somewhat more difficult.

\bigskip

The following two theorems, for the commutative and the free associative cases, respectively, constitute the foundation of the approximation technique.

\begin{theorem}   \label{ThMainTechn}
Let $\varphi\in\Aut_{\Ind}(\Aut_0(K[x_1,\dots,x_n]))$ and let $H_N\subset \Aut_0(K[x_1,\dots,x_n])$ be the subgroup of
automorphisms which are identity modulo the ideal
 $(x_1,\dots,x_n)^N$ ($N>1$). Then $\varphi(H_N)\subseteq H_N$.

\end{theorem}

\begin{theorem}   \label{ThMainTechnFreeAss}
Let $\varphi\in\Aut_{\Ind}(\Aut_0(K\langle x_1,\dots,x_n\rangle ))$ and let $H_N$ be again the subgroup of
automorphisms which are identity modulo the ideal
 $(x_1,\dots,x_n)^N$. Then $\varphi(H_N)\subseteq H_N$.
\end{theorem}

\begin{corollary}
In both commutative and free associative cases under the assumptions above one has $\varphi= \Id$.
\end{corollary}

{\bf Proof.} Every automorphism can be approximated via the tame ones,
i.e. for any $\psi$ and any $N$ there exists a tame  automorphism
$\psi'_N$ such that $\psi\psi'_N{}^{-1}\in H_N$.


The main point therefore is why $\varphi(H_N)\subseteq H_N$ whenever $\varphi$ is and $\Ind$-automorphism.

\medskip

{\bf Proof of Theorem \ref{ThMainTechn}.}

The method of proof is based upon the following useful fact from algebraic geometry:

\begin{lemma} \label{oneparamlem} Let

$$
\varphi: X\rightarrow Y
$$
be a morphism of affine varieties, and let $A(t)\subset X$ be a curve (or rather, a one-parameter family of points) in $X$. Suppose that $A(t)$ does not tend to infinity as $t\rightarrow 0$. Then the image $\varphi A(t)$ under $\varphi$ also does not tend to infinity as $t\rightarrow 0$.
\end{lemma}

The proof is straightforward and is left to the reader.

We now put the above fact to use. For $t>0$ let
$$
\hat{A}(t):\mathbb{A}^n_{K}\rightarrow \mathbb{A}^n_{K}
$$
be a one-parameter family of invertible linear transformations of the affine space preserving the origin. To that corresponds a curve
$A(t)\subset \Aut_0(K[x_1,\dots,x_n])$ of polynomial automorphisms whose points are linear substitutions. Suppose that, as $t$ tends to zero, the $i$-th eigenvalue of $A(t)$ also tends to zero as $t^{k_i}$, $k_i\in\mathbb{N}$. Such a family will always exist.

Suppose now that the degrees $\lbrace k_i,\;i=1,\ldots n\rbrace$ of singularity of eigenvalues at zero are such that for every pair $(i,j)$, if $k_i\neq k_j$, then there exists a positive integer $m$ such that
$$
\text{either\;\;} k_im\leq k_j\;\;\text{or\;\;}k_jm\leq k_i.
$$

The largest such $m$ we will call the order of $A(t)$ at $t=0$. As $k_i$ are all set to be positive integer, the order equals $\frac{k_{\text{max}}}{k_{\text{min}}}$.

Let $M\in \Aut_0(K[x_1,\dots,x_n])$ be a polynomial automorphism.
\begin{lemma}    \label{Lm2} The curve $A(t)MA(t)^{-1}$ has no singularity at zero for any $A(t)$ of order $\leq N$ if and only if $M\in \hat{H}_N$, where $\hat{H}_N$ is the subgroup of automorphisms which are homothety modulo the
augmentation ideal.
\end{lemma}






{\bf Proof.} The  `If' part is elementary, for if $M\in \hat{H}_N$, the action of $A(t)MA(t)^{-1}$ upon any generator $x_i$ (with $i$ fixed)\footnote{Without loss of generality we may assume that the coordinate functions $x_i$ correspond to the principal axes of $\hat{A}(t)$.} is given by
\begin{eqnarray*}
A(t)MA(t)^{-1}(x_i)=\lambda x_i+t^{-k_i}\sum_{l_1+\cdots +l_n=N}a_{l_1\ldots l_n}t^{k_1l_1+\cdots+k_nl_n}x_1^{l_1}\cdots x_n^{l_n}+\\
+S_i(t,x_1,\ldots,x_n),
\end{eqnarray*}
where $\lambda$ is the homothety ratio of (the linear part of) $M$ and $S_i$ is polynomial in $x_1,\ldots,x_n$ of total degree greater than $N$. Now, for any choice of $l_1,\ldots,l_n$ in the sum, the expression
$$
k_1l_1+\cdots+k_nl_n-k_i\geq k_{\text{min}}\sum l_j-k_i=k_{\text{min}}N-k_i\geq 0
$$
for every $i$, so whenever $t$ goes to zero, the coefficient will not blow up to infinity. Obviously the same argument applies to higher-degree monomials within $S_i$.

\medskip

The other direction is slightly less elementary; assuming that $M\notin \hat{H}_N$, we need to show that there is a curve $A(t)$ such that conjugation of $M$ by it produces a singularity at zero. We distinguish between two cases.

{\bf Case 1.} The linear part $\bar{M}$ of $M$ is not a scalar
matrix. Then -- after a suitable basis change (see the footnote) - it is not a diagonal matrix and
has a non-zero entry in the position $(i,j)$. Consider a
diagonal matrix $A(t)=D(t)$ such that on all positions on the main
diagonal except $j$-th it has $t^{k_i}$ and on $j$-th position it has
$t^{k_j}$. Then $D(t)\bar{M}D^{-1}(t)$ has $(i,j)$ entry with the
coefficient $t^{k_i-k_j}$  and if $k_j>k_i$ it has a
singularity at $t=0$.

Let also $k_i<2k_j$. Then the non-linear part of $M$ does not produce
singularities and cannot compensate the singularity of the linear part.

{\bf Case 2.} The linear part $\bar{M}$ of $M$ is  a scalar
matrix. Then conjugation cannot produce
singularities in the linear part and we as before are interested in the smallest non-linear
term. Let $M\in H_N\backslash H_{N+1}$. Performing a basis change if necessary,
we may assume that
$$
\varphi(x_1)=\lambda x_1+\delta
x_2^N+S,
$$
where $S$ is a sum of monomials of degree $\ge N$
with coefficients in $K$.

Let $A(t)=D(t)$ be a diagonal matrix of the form $(t^{k_1},
t^{k_2},t^{k_1},\dots,t^{k_1})$ and let $(N+1)\cdot k_2>k_1>N\cdot
k_2$. Then in $A^{-1}MA$ the term $\delta x_2^N$ will be transformed
into $\delta x_2^N t^{Nk_2-k_1}$, and all other terms are multiplied by
$t^{lk_2+sk_1-k_1}$ with $(l,s)\ne (1,0)$ and $ l, s>0$. In this
case $lk_2+sk_1-k_1>0$ and we are done with the proof of  Lemma
\ref{Lm2}.

\medskip

The next lemma is proved by direct computation. Recall that for $m>1$, the group $G_m$ is defined as the group of all tame automorphisms preserving the $m$-th power of the augmentation ideal.

\begin{lemma}     \label{Lm3}

\

a) $[G_m,G_m]\subset H_m$, $m>2$. There exist elements \\
$\varphi\in
H_{m+k-1}\backslash H_{m+k},\;\;\psi_1\in G_k,\;\;\psi_2\in G_m$, such that
$\varphi=[\psi_1,\psi_2]$.

b) $[H_m,H_k]\subset H_{m+k-1}$.

c) Let $\varphi\in G_m\backslash H_{m}$, $\psi\in H_k\backslash
H_{k+1}$, $k>m$. Then $[\varphi,\psi]\in H_k\backslash H_{k+1}$.

\end{lemma}

{\bf Proof.} a) Consider elementary automorphisms
\begin{gather*}
\psi_1: x_1\mapsto x_1+x_2^k,\;\;
 x_2\mapsto x_2,\;\;
 x_i\mapsto x_i,\;i> 2;
\end{gather*}
\begin{gather*}
\psi_2: x_1\mapsto x_1,\;\;
x_2\mapsto x_2+x_1^m, \;\;
x_i\mapsto x_i,\;i>2.
\end{gather*}

Set $\varphi=[\psi_1,\psi_2]=\psi_1^{-1}\psi_2^{-1}\psi_1\psi_2$. \\
Then
$$
\varphi: x_1\mapsto x_1-x_2^k+(x_2-(x_1-x_2^k)^m)^k,$$ $$x_2\mapsto
x_2-(x_1-x_2^k)^m+(x_1-x_2^k+(x_2-(x_1-x_2^k)^m)^k)^m, \;\;
x_i\mapsto x_i,\;i>2.
$$
It is easy to see that if either $k$ or
$m$ is relatively prime with $\Ch(K)$, then not all terms of degree
$k+m-1$ vanish. Thus $\varphi\in H_{m+k-1}\backslash
H_{m+k}$.

Now suppose that $\Ch(K)\nmid m$, then obviously $m-1$ is
relatively prime with $\Ch(K)$. Consider the mappings
$$
\psi_1: x_1\mapsto x_1+x_2^k,\;\; x_2\mapsto x_2,\;\;x_i\mapsto x_i,\;i>2;
$$
$$
\psi_2: x_1\mapsto x_1,\;\; x_2\mapsto x_2+x_1^{m-1}x_3,\;\;x_i\mapsto x_i,\;i>2.
$$
Set again
$\varphi'=[\psi_1,\psi_2]=\psi_1^{-1}\psi_2^{-1}\psi_1\psi_2$.
Then $\varphi'$ acts as
\begin{gather*}
x_1\mapsto x_1-x_2^k+(x_2-(x_1-x_2^k)^{m-1}x_3)^k= \\
=x_1-k(x_1-x_2^k)^{m-1}x_2^{k-1}x_3+S, \\
x_2\mapsto x_2-(x_1-x_2^k)^{m-1}x_3+(x_1-x_2^k+(x_2-(x_1-x_2^k)^{m-1}x_3)^k)^{m-1}x_3, \\
x_i\mapsto x_i,\;\; i>2;
\end{gather*}
here $S$ stands for a sum of terms of degree $\ge m+k$. Again we see that $\varphi\in
H_{m+k-1}\backslash H_{m+k}$.

b) Let $$\psi_1: x_i\mapsto x_i+f_i;\ \psi_2: x_i\mapsto x_i+g_i,$$
for $i=1,\dots,n;$ here $f_i$ and $g_i$ do not have monomials of degree less than or equal to $m$ and $k$, respectively. Then, modulo terms
of degree $\ge m+k$, we have $\psi_1\psi_2: x_i\mapsto
x_i+f_i+g_i+\frac{\partial f_i}{\partial x_j} g_j$, so that modulo
terms of degree $\ge m+k-1$ we get $\psi_1\psi_2: x_i\mapsto
x_i+f_i+g_i$ and $\psi_2\psi_1: x_i\to x_i+f_i+g_i$. Therefore $[\psi_1,\psi_2]\in H_{m+k-1}$.

c) If $\varphi(I^m)\subseteq I^m$ and
$$
\psi: (x_1,\ldots,x_n)\mapsto (x_1+g_1,\ldots,x_n+g_n)
$$
is such that for some $i_0$ the polynomial $g_{i_0}$ contains a monomial of total degree $k$ (and all $g_i$ do not contain monomials of total degree less than $k$), then, by evaluating the composition of automorphisms directly, one sees that the commutator is given by
$$
[\varphi,\psi]: (x_1,\ldots,x_n)\mapsto (x_1+g_1+S_1,\ldots,x_n+g_n+S_n)
$$
with $S_i$ containing no monomials of total degree $<k+1$. Then the image of $x_{i_0}$ is $x_{i_0}$ modulo polynomial of height $k$.

\begin{corollary}      \label{Cofinal}
Let $\Psi\in\Aut_{\Ind}(\NAut(K[x_1,\dots,x_n]))$. Then\\
$\Psi(G_n)=G_n$, $\Psi(H_n)=H_n$.
\end{corollary}

Corollary \ref{Cofinal} together with Proposition \ref{PrRips} of the next section imply Theorem
\ref{ThMainTechn}, for every nice automorphism, by definition, can be
approximated by tame ones. Note that in characteristic zero
every automorphism is nice (Anick's theorem).

\subsection{Lifting of automorphism groups}

\subsubsection{Lifting of automorphisms from
$\Aut(K[x_1,\dots,x_n])$ to
$\Aut(K\langle x_1,\dots,x_n\rangle)$}
\begin{definition} \label{defBialBir}
In the sequel, we call an action of the $n$-dimensional torus ${\mathbb T}^n$ on ${\mathbb
K}\langle x_1,\dots,x_n\rangle$ (the number of generators coincides with the dimension of the torus) \textbf{linearizable} if it is conjugate to the standard diagonal action given by
$$
(\lambda_1,\ldots,\lambda_n)\;:\;(x_1,\ldots,x_n)\mapsto (\lambda_1x_1,\ldots,\lambda_nx_n).
$$
\end{definition}
The following result is a direct free associative analogue of a well-known theorem of Bia\l{}ynicki-Birula \cite{Bialynicki-Birula2,Bialynicki-Birula1}. We will make frequent reference of the classical (commutative) case as well, which appears as Theorem \ref{ThBialBir} in the text.

\begin{theorem}  \label{ThBialBirFree}
Any effective action of the $n$-torus on ${\mathbb
K}\langle x_1,\dots,x_n\rangle$ is linearizable.
\end{theorem}

The proof is somewhat similar to that of Theorem \ref{ThBialBir}, with a few modifications. We will address this issue in the upcoming paper \cite{KBBialBir}.

\medskip

As a corollary of the above theorem,  we get

\begin{proposition}
Let $T^n$ denote the standard torus action on
$K[x_1,\ldots,x_n]$. Let $\widehat{T}^n$ denote its
lifting to an action on the free associative algebra $K\langle x_1,\ldots,x_n\rangle$. Then
$\widehat{T}^n$ is also given by the standard torus action.
\end{proposition}

{\bf Proof.} Consider the roots $\widehat{x_i}$ of this action.
They are liftings of the coordinates $x_i$. We have to prove that
they generate the whole associative algebra.

Due to the reducibility of this action, all elements are product
of eigenvalues of this action. Hence it is enough to prove that
eigenvalues of this action can be presented as a linear combination
of this action. This can be done along the lines of Bia\l{}ynicki-Birula
\cite{Bialynicki-Birula1}. Note that all propositions of the previous
section hold for the free associative algebra. Proof of Theorem
\ref{ThMainTechnFreeAss} is  similar. Hence we have the
following

 \begin{theorem}        \label{ThAutAutFree}
Any $\Ind$-scheme automorphism $\varphi$ of
$\Aut(K\langle
x_1,\dots,x_n\rangle)$ for $n\ge 3$ is inner, i.e. is a
conjugation by some automorphism.
\end{theorem}

We therefore see that the group lifting (in the sense of isomorphism induced by
the natural abelianization) implies the analogue of Theorem
\ref{ThMainTechn}.

This also implies that any automorphism group lifting, if exists, satisfies the
approximation properties.

\begin{proposition}    \label{PrLiftStrct}
Suppose
$$
\Psi: \Aut(K[x_1,\dots,x_n])\rightarrow \Aut(K\langle
z_1,\dots,z_n\rangle)
$$
is a group
homomorphism such that its composition with the natural map $\Aut(K\langle
z_1,\dots,z_n\rangle)\rightarrow \Aut(K[x_1,\dots,x_n])$ (induced by the projection $K\langle
z_1,\dots,z_n\rangle\rightarrow K[x_1,\dots,x_n]$) is
the identity map. Then

\begin{enumerate}
    \item After a coordinate change $\Psi$ provides a
    correspondence between the standard torus actions $x_i\mapsto
    \lambda_ix_i$ and  $z_i\mapsto
    \lambda_iz_i$.
    \item Images of elementary automorphisms
    $$x_j\mapsto x_j,\; j\ne i,\;\; x_i\mapsto x_i+f(x_1,\dots,\widehat{x_i},\dots,x_n)$$
    are elementary
    automorphisms of the form
        $$z_j\mapsto z_j,\; j\ne i,\;\; z_i\mapsto z_i+f(z_1,\dots,\widehat{z_i},\dots,z_n).$$
    (Hence image of tame automorphism is tame
    automorphism).
    \item $\psi(H_n)=G_n$. Hence $\psi$ induces a map between the
    completion of the groups of $\Aut(K[x_1,\dots,x_n])$ and $\Aut(K\langle
z_1,\dots,z_n\rangle)$ with respect to the
    augmentation subgroup structure.
\end{enumerate}
\end{proposition}

\

\noindent {\bf Proof of Theorem \ref{ThGroupLifting}}

\noindent Any automorphism (including wild automorphisms such as the Nagata example) can
be approximated by a product of elementary automorphisms with
respect to augmentation topology. In the case of the Nagata
automorphism corresponding to
$$\Aut(K\langle
x_1,\dots,x_n\rangle),$$ all such elementary automorphisms fix all
coordinates except  $x_1$ and $x_2$. Because of (2) and (3) of
Proposition \ref{PrLiftStrct}, the lifted automorphism would be an
automorphism induced by an automorphism of $K\langle
x_1,x_2,x_3\rangle$ fixing $x_3$. However, it is impossible to
lift the Nagata automorphism to such an automorphism due to the
main result of \cite{BelovYuLifting}. Therefore, Theorem
\ref{ThGroupLifting} is proved.



\section{Automorphisms of the polynomial algebra and the approach of Bodnarchuk--Rips} \label{ScAutTameCoommN2}
 Let
$\Psi\in\Aut(\Aut(K[x_1,\dots,x_n]))$ (resp.
$\Aut(\TAut(K[x_1,\dots,x_n]))$,\\
$\Aut(\TAut_0(K[x_1,\dots,x_n]))$,
$\Aut(\Aut_0(K[x_1,\dots,x_n]))$).

\subsection{Reduction to the case when $\Psi$ is identical on
$\SL_n$} We follow \cite{KraftStampfli} and
\cite{Bodnarchuk} using the classical theorem of Bia\l{}ynicki-Birula
\cite{Bialynicki-Birula2,Bialynicki-Birula1}:

\begin{theorem}[Bia\l{}ynicki-Birula]  \label{ThBialBir}
Any effective action of torus ${\mathbb T}^n$ on ${\mathbb C}^n$
is linearizable (recall the definition \ref{defBialBir}).
\end{theorem}

{\small

{\bf Remark.} An effective action of ${\mathbb T}^{n-1}$ on
${\mathbb C}^n$ is linearizable
\cite{Bialynicki-Birula1,Bialynicki-Birula2}. There is a
conjecture whether any action of ${\mathbb T}^{n-2}$ on ${\mathbb
C}^n$ is linearizable, established for $n=3$. For codimension
 $>2$, there are positive-characteristic counterexamples \cite{Asanuma}.

\medskip

{\bf Remark.} Kraft and Stampfli \cite{KraftStampfli} proved
(by considering periodic elements in $\mathbb T$) that an effective action $T$ has the following property: if $\Psi\in\Aut(\Aut)$
is a group automorphism, then the image of $T$ (as a subgroup of $\Aut$) under $\Psi$ is an algebraic group. In fact their proof
is also applicable for the free associative algebra case.
We are going to use this result.

\medskip
}

Returning to the case of automorphisms $\varphi\in\Aut_{\Ind}\Aut$ preserving the $\Ind$-group structure, consider now
the standard action $x_i\mapsto\lambda_ix_i$ of the $n$-dimensional torus $\mathbb{T}\leftrightarrow T^n\subset \Aut(\mathbb{C}[x_1,\ldots,x_n])$ on the affine space $\mathbb{C}^n$.
Let $H$ be the image of $T^n$ under $\varphi$. Then
by Theorem \ref{ThBialBir} $H$ is conjugate to the standard
torus $T^n$ via some automorphism $\psi$. Composing $\varphi$ with this
conjugation, we come to the case when $\varphi$
is the identity on the maximal torus.  Then we have the following

\begin{corollary}
Without loss of generality, it is enough to prove Theorem
\ref{ThAutAut} for the case when $\varphi|_{\mathbb T}=\Id$.
\end{corollary}

Now we are in the situation when $\varphi$ preserves all linear
mappings $x_i\mapsto \lambda_i x_i$.  We have to prove that it is
the identity.

\begin{proposition}[E. Rips, private communication] \label{PrRips}
Let $n>2$ and suppose $\varphi$ preserves the standard torus action
on the commutative polynomial algebra. Then $\varphi$ preserves all elementary
transformations.
\end{proposition}

\begin{corollary}
Let $\varphi$ satisfy the conditions of Proposition
\ref{PrRips}. Then $\varphi$ preserves all tame automorphisms.
\end{corollary}

{\bf Proof of Proposition \ref{PrRips}.} We state a few elementary lemmas.

\begin{lemma}  \label{LmTorComp}
Consider the diagonal action $T^1\subset T^n$ given by automorphisms: $\alpha:
x_i\mapsto \alpha_i x_i$, $\beta: x_i\mapsto \beta_i x_i$. Let $\psi:
x_i\mapsto \sum_{i,J} a_{iJ}x^J,\; i=1,\dots,n$, where $J=(j_1,\dots,j_n)$ is the
multi-index, $x^J=x^{j_1}\cdots x^{j_n}$. Then

$$\alpha\circ\psi\circ\beta: x_i\mapsto \sum_{i,J} \alpha_i
a_{iJ}x^J\beta^J,$$

In particular,
$$\alpha\circ\psi\circ\alpha^{-1}: x_i\mapsto \sum_{i,J} \alpha_i
a_{iJ}x^J\alpha^{-J}.$$
\end{lemma}

Applying Lemma \ref{LmTorComp} and comparing the coefficients we
get the following

\begin{lemma}  \label{LmLindiag}
Consider the diagonal $T^1$ action: $x_i\mapsto \lambda x_i$.
Then the set of automorphisms commuting with this action is
exactly the set of linear automorphisms.
\end{lemma}

Similarly (using Lemma \ref{LmTorComp}) we obtain Lemmas
\ref{LmMult1}, \ref{LmMult1n}, \ref{LmMult2}:

\begin{lemma}     \label{LmMult1}
a) Consider the following $ T^2$ action: $$x_1\mapsto
\lambda\delta x_1,\; x_2\mapsto \lambda x_2,\; x_3\mapsto \delta x_3,\;
x_i\mapsto \lambda x_i,\; i>3.$$ Then the set $S$ of automorphisms commuting with
this action is generated by the following automorphisms: $$x_1\mapsto
x_1+\beta x_2x_3,\; x_i\mapsto \varepsilon_i x_i,\; i>1,\;
(\beta,\varepsilon_i\in K).$$

b) Consider the following $T^{n-1}$ action:
$$x_1\mapsto
\lambda^I x_1,\; x_j\mapsto \lambda_j x_j,\; j>1\;(\lambda^I=\lambda_2^{i_2}\cdots\lambda_n^{i_n}).$$ Then the set $S$ of
automorphisms commuting with this action is generated by the following
automorphisms: $$x_1\mapsto x_1+\beta \prod_{j=2}^n x_j^{i_j},\
(\beta \in K).$$
\end{lemma}

{\bf Remark.} A similar statement for the free associative case is
true, but one has to consider the set $\hat{S}$ of automorphisms
$x_1\mapsto x_1+h,\; x_i\mapsto \varepsilon_i x_i,\; i>1$, ($\varepsilon\in
K$, and the polynomial $h\in K\langle x_2,\dots,x_n\rangle$ has
total degree $J$ - in the free associative case it is not just monomial
anymore).

\begin{corollary}      \label{CoActionToronElem}
Let $\varphi\in\Aut(\TAut(K[x_1,\dots,x_n]))$  stabilizing all
elements from $\mathbb T$. Then $\varphi(S)=S$.
\end{corollary}

\begin{lemma}     \label{LmMult1n}
Consider the following $T^1$ action:
$$
x_1\mapsto \lambda^2 x_1,\;
x_i\mapsto \lambda x_i,\; i>1.
$$
Then the set $S$ of
automorphisms commuting with this action is generated by the following
automorphisms:
$$
x_1\mapsto x_1+\beta x_2^2,\; x_i\mapsto  \lambda_i x_i,\;
i>2,\; (\beta,\lambda_i\in K).
$$
\end{lemma}

\begin{lemma} \label{LmMult2}
Consider the set $S$ defined in the previous lemma. Then
$[S,S]=\{uvu^{-1}v^{-1}\}$ consists of the following automorphisms
$$x_1\mapsto x_1+\beta x_2x_3,\; x_2\mapsto x_2,\; x_3\mapsto  x_3,
\;(\beta\in K).$$
\end{lemma}

\begin{lemma}  \label{LmMult3}
Let $n\ge 3$. Consider the following set of automorphisms $$\psi_i:
x_i\mapsto x_i+\beta_ix_{i+1}x_{i+2},\; \beta_i\ne 0,\; x_k\mapsto x_k,\; k\ne i$$
for $i=1,\dots, n-1$. (Numeration is cyclic, so for example
$x_{n+1}=x_1$). Let $\beta_i\ne 0$ for all $i$. Then all of
$\psi_i$ can be simultaneously conjugated by a torus action to $$\psi_i':
x_i\mapsto x_i+x_{i+1}x_{i+2},\; x_k\mapsto x_k,\; k\ne i$$ for $i=1,\dots, n$ in
a unique way.
\end{lemma}

{\bf Proof.} Let $\alpha: x_i\mapsto \alpha_i x_i$. Then by Lemma
\ref{LmTorComp} we obtain
$$\alpha\circ\psi_i\circ\alpha^{-1}: x_i\mapsto x_i+\beta_i
x_{i+1}x_{i+2}\alpha_{i+1}^{-1}\alpha_{i+2}^{-1}\alpha_i$$ and
$$\alpha\circ\psi_i\circ\alpha^{-1}: x_k\mapsto x_k$$
for $k\ne i$.

Comparing the coefficients of the quadratic terms, we see that it is
sufficient to solve the system:
$$\beta_i\alpha_{i+1}^{-1}\alpha_{i+2}^{-1}\alpha_i=1,\; i=1,\dots,n-1.$$
As $\beta_i\ne 0$ for all $i$, this system has a unique
solution.

{\bf Remark.} In the free associative algebra case, instead of $\beta x_2x_3$ one
has to consider $\beta x_2x_3+\gamma x_3x_2$.

\subsection{The  lemma of Rips}

\begin{lemma}[E. Rips]  \label{LmRips}
Let $\Ch(K)\ne 2$, $|K|=\infty$. Linear transformations and
$\psi_i'$ defined in Lemma \ref{LmMult3} generate the whole tame automorphism group of $K[x_1,\ldots,x_n]$.
\end{lemma}

Proposition \ref{PrRips}  follows from Lemmas \ref{LmLindiag},
\ref{LmMult1}, \ref{LmMult1n}, \ref{LmMult2}, \ref{LmMult3},
\ref{LmRips}. Note that we have proved an analogue of Theorem
\ref{ThAutAut} for tame automorphisms.

\medskip
{\bf Proof of Lemma \ref{LmRips}.} Let $G$ be the group generated
by elementary transformations as in Lemma \ref{LmMult3}. We
have to prove that is isomorphic to the tame automorphism subgroup fixing
the augmentation ideal. We are going to need some preliminaries.

\begin{lemma}  \label{LmR1}
Linear transformations of $K^3$ and $$\psi: x\mapsto x,\; y\mapsto y,\; z\mapsto z+xy$$
generate all mappings of the form $$\phi_m^b(x,y,z):
x\mapsto x,\; y\mapsto y,\; z\mapsto z+bx^m,\;\; b\in K.$$
\end{lemma}

{\bf Proof of Lemma \ref{LmR1}.} We proceed by induction.
Suppose we have an automorphism $$\phi^b_{m-1}(x,y,z): x\mapsto x,\; y\mapsto
y,\; z\mapsto z+bx^{m-1}.$$  Conjugating by the linear transformation
($z\mapsto y,\; y\mapsto z,\;x\mapsto x$), we obtain the automorphism
$$\phi_{m-1}^b(x,z,y): x\mapsto x,\; y\mapsto y+bx^{m-1},\; z\mapsto z.$$
Composing this on the right by $\psi$, we get the automorphism
$$\varphi(x,y,z): x\mapsto x,\; y\mapsto y+bx^{m-1},\; z\mapsto z+yx+x^m.$$ Note
that
$$\phi_{m-1}(x,y,z)^{-1}\circ\varphi(x,y,z): x\mapsto x,\; y\mapsto y,\;
z\mapsto z+xy+bx^m.$$
Now we see that
$$\psi^{-1}\phi_{m-1}(x,y,z)^{-1}\circ\varphi(x,y,z)=\phi^b_m$$ and
the lemma is proved.

\begin{corollary}  \label{CoLmR1}
Let $\Ch(K)\nmid n$ (in particular, $\Ch(K)\ne 0$) and
$|K|=\infty$. Then $G$ contains all the transformations $$z\mapsto
z+bx^ky^l,\; y\mapsto y,\; x\mapsto x$$ such that $k+l=n$.
\end{corollary}

{\bf Proof.} For any invertible linear transformation
$$\varphi:
x\mapsto a_{11}x+a_{12}y,\; y\mapsto a_{21}x+a_{22}y,\; z\mapsto z; a_{ij}\in K$$
we have $$\varphi^{-1}\phi^b_{m}\varphi: x\mapsto x,\; y\mapsto y,\; z\mapsto
z+b(a_{11}x+a_{12}y)^m.$$ Note that sums of such expressions
contain all the terms of the form $bx^ky^l$. The corollary is proved.

\subsection{Generators of the tame automorphism group}

\begin{theorem}  \label{ThR2}
If $\Ch(K)\ne 2$ and $|K|=\infty$, then linear transformations and
$$\psi: x\mapsto x,\; y\mapsto y,\; z\mapsto z+xy$$ generate all mappings of the
form $$\alpha_m^b(x,y,z): x\mapsto x,\; y\mapsto y,\; z\mapsto z+byx^m,\;\;
b\in K.$$
\end{theorem}

{\bf Proof of theorem \ref{ThR2}.} Observe that
$$\alpha=\beta\circ\phi^b_m(x,z,y): x\mapsto x+by^m,\; y\mapsto y+x+by^m,\;
z\mapsto z,$$ where $\beta: x\mapsto x,\; y\mapsto x+y,\;z\mapsto z$. Then
$$\gamma=\alpha^{-1}\psi\alpha: x\mapsto x,\; y\mapsto y,\; z\mapsto
z+xy+2bxy^m+by^{2m}.$$ Composing with $\psi^{-1}$ and
$\phi_{2m}^{2b}$ we get the desired $$\alpha_m^{2b}(x,y,z): x\mapsto x,\; y\mapsto
y,\; z\mapsto z+2byx^m,\;\; b\in K.$$

\begin{corollary}  \label{CoThR2}
Let $\Ch(K)\nmid n$ and $|K|=\infty$. Then $G$
contains all transformations of the form $$z\mapsto z+bx^ky^l,\; y\mapsto y,\; x\mapsto x$$
such that $k=n+1$.
\end{corollary}

The {\bf proof} is similar to the proof of Corollary
\ref{CoLmR1}. Note that either $n$ or $n+1$ is not a multiple of $\Ch(K)$ so
we have

\begin{lemma}  \label{LmR2n}
If $\Ch(K)\ne 2$ then linear transformations and $$\psi: x\mapsto x,\;
y\mapsto y,\; z\mapsto z+xy$$ generate all mappings of the form
$$\alpha_P: x\mapsto x,\; y\mapsto y,\; z\mapsto z+P(x,y),\;\; P(x,y)\in K[x,y].$$
\end{lemma}

We have proved Lemma \ref{LmRips} for the three variable case. In order
to treat the case $n\ge 4$ we need one more lemma.

\begin{lemma}  \label{LmR3Lin}
Let $M(\vec x)=a\prod x_i^{k_i},\;$ $a\in K,\; |K|=\infty$,
$\Ch(K)\nmid k_i$ for at least one of $k_i$'s. Consider the linear
transformations denoted by $$f: x_i\mapsto y_i=\sum a_{ij}x_j,\; \det(a_{ij})\ne
0$$ and monomials $M_f=M(\vec y)$. Then the linear span of $M_f$ for
different $f$'s contains all homogenous polynomials of degree
$k=\sum k_i$ in $K[x_1,\dots,x_n]$.
\end{lemma}

{\bf Proof.} It is a direct consequence of the following fact.
Let $S$ be a homogenous subspace of $K[x_1,\dots,x_n]$ invariant  with
respect to $GL_n$ of degree $m$. Then $S=S_{m/p^k}^{p^k}$,
$p=\Ch(K)$,\ $S_l$ is the space of all polynomials of degree $l$.

Lemma \ref{LmRips} follows from Lemma \ref{LmR3Lin} in a similar way as in the
proofs of Corollaries \ref{CoLmR1} and \ref{CoThR2}.

\

\subsection{$\Aut(\TAut)$ for general case}

 Now we consider the case when  $\Ch(K)$ is arbitrary, i.e. the remaining case
 $\Ch(K)=2$. Still $|K|=\infty$.
Although we are unable to prove the analogue of Proposition
\ref{PrRips}, we can still play on the relations.

Let $$M=a\prod_{i=1}^{n-1} x_i^{k_i}$$ be a monomial, $a\in K$. For
polynomial $P(x,y)\in K[x,y]$ we define the elementary automorphism
$$\psi_P: x_i\mapsto x_i,\; i=1,\dots,n-1,\; x_n\mapsto
x_n+P(x_1,\dots,x_{n-1}).$$ We have $P=\sum M_j$ and $\psi_P$
naturally decomposes as a product of commuting $\psi_{M_j}$. Let
$\Psi\in\Aut(\TAut(K[x,y,z]))$ stabilizing linear mappings and
$\phi$ (Automorphism $\phi$ defined in Lemma \ref{LmR1}). Then
according to the corollary \ref{CoActionToronElem}
$\Psi(\psi_P)=\prod \Psi(\psi_{M_j})$. If $M=ax^n$ then due to
Lemma \ref{LmR1}
$$\Psi(\psi_{M})=\psi_{M}.$$ We have to prove the same for other
type of monomials:

\begin{lemma}  \label{Le2nTestFreeAss}
Let $M$ be a monomial. Then
$$\Psi(\psi_{M})=\psi_{M}.$$
\end{lemma}

{\bf Proof.} Let $M=a\prod_{i=1}^{n-1} x_i^{k_i}$. Consider the
automorphism $$\alpha: x_i\mapsto x_i+x_1,\; i=2,\dots,n-1;\; x_1\mapsto
x_1,\;x_n\mapsto x_n.$$ Then $$\alpha^{-1}\psi_M\alpha
=\psi_{x_1^{k_1}\prod_{i=2}^{n-1}(x_i+x_1)^{k_i}}=\psi_Q\psi_{ax_1^{\sum_{i=2}^{n-1}
k_i}}.$$

Here the polynomial
$$Q=x_1^{k_1}\left(\prod_{i=2}^{n-1}(x_i+x_1)^{k_i}-ax_1^{\sum
k_i}\right).$$

It has the following form
$$Q=\sum_{i=2}^{n-1} N_i,$$
where $N_i$ are monomials such that none of them is proportional to a
power of $x_1$.

According to Corollary \ref{CoActionToronElem},
$\Psi(\psi_M)=\psi_{bM}$ for some $b\in K$. We need only to prove
that $b=1$. Suppose the contrary, $b\ne 1$. Then
\begin{eqnarray*}
\Psi(\alpha^{-1}\psi_M\alpha)=\left(\prod_{[N_i,x_1]\ne 0}\Psi(\psi_{N_i})\right)\circ\Psi(\psi_{ax_1^{\sum_{i=2}^{n-1}
k_i}})=\\
\left(\prod_{[N_i,x_1]\ne
0}\psi_{b_iN_i}\right)\circ\psi_{ax_1^{\sum_{i=2}^{n-1} k_i}}
\end{eqnarray*}
for some $b_i\in K$.

On the other hand
$$\Psi(\alpha^{-1}\psi_M\alpha)=\alpha^{-1}\Psi(\psi_M)\alpha=\alpha^{-1}\psi_{bM}\alpha =
\left(\prod_{[N_i,x_1]\ne 0}
\psi_{bN_i}\right)\circ\psi_{ax_1^{\sum_{i=2}^{n-1} k_i}}$$

Comparing the factors $\psi_{ax_1^{\sum_{i=2}^{n-1} k_i}}$ and
$\psi_{ax_1^{\sum_{i=2}^{n-1} k_i}}$ in the last two products we
get $b=1$. Lemma \ref{Le2nTestFreeAss} and hence Proposition
\ref{PrRips} are proved.

\section{The approach of Bodnarchuk--Rips to automorphisms of $\TAut(K\langle
x_1,\dots,x_n\rangle)$ ($n>2$)}

Now consider the free associative case.
We treat the case $n>3$ on group-theoretic
level and the case $n=3$ on $\Ind$-scheme level. Note that if $n=2$
then $\Aut_0(K[x,y])=\TAut_0(K[x,y])\simeq \TAut_0(K\langle
x,y\rangle)=\Aut_0(K\langle x,y\rangle)$ and description of
automorphism group of such objects is known due to J. D\'{e}serti.

\subsection{The automorphisms of the tame automorphism \\
group of
$K\langle x_1,\dots,x_n\rangle$, $n\ge 4$}

\begin{proposition}[E. Rips, private communication] \label{PrRipsFreeAss}
 Let $n>3$ and let $\varphi$ preserve the standard torus action
on the free associative  algebra $K\langle x_1,\dots,x_n\rangle$. Then $\varphi$ preserves all
elementary transformations.
\end{proposition}

\begin{corollary}
Let $\varphi$ satisfy the conditions of the proposition
\ref{PrRipsFreeAss}. Then $\varphi$ preserves all tame
automorphisms.
\end{corollary}

For free associative algebras, we note that any automorphism
preserving the torus action preserves also the symmetric $$x_1\mapsto
x_1+\beta(x_2x_3+x_3x_2),\; x_i\mapsto x_i,\; i>1$$ and the skew symmetric
$$x_1\mapsto x_1+\beta(x_2x_3-x_3x_2),\; x_i\mapsto x_i,\; i>1$$ elementary
automorphisms. The first property follows from Lemma \ref{LmMult1n}.
The second one follows from the fact that skew symmetric automorphisms
commute with automorphisms of the following type $$x_2\mapsto x_2+x_3^2,\;
x_i\mapsto x_i,\; i\ne 2$$ and this property distinguishes them from
elementary automorphisms of the form $$x_1\mapsto x_1+\beta
x_2x_3+\gamma x_3x_2,\; x_i\mapsto x_i,\; i>1.$$

Theorem \ref{ThAutAutFreeAssoc} follows from the fact that the
forms $\beta x_2x_3+\gamma x_3x_2$ corresponding to general bilinear multiplication $$*_{\beta,\gamma}:(x_2,\;x_3)\to \beta x_2x_3+\gamma x_3x_2$$
lead to associative multiplication if and only if $\beta=0$ or $\gamma=0$;
the approximation also applies (see section \ref{SbScAprox}).

\medskip
Suppose at first that $n=4$ and we are dealing with $K\langle x,y,z,t\rangle$.

\begin{proposition}       \label{PrGTame2}
The group $G$ containing all linear transformations and mappings
$$x\mapsto x,\; y\mapsto y,\; z\mapsto z+xy,\; t\mapsto t$$ contains also all
transformations of the form $$x\mapsto x,\; y\mapsto y,\; z\mapsto z+P(x,y),\; t\mapsto t.$$
\end{proposition}

{\bf Proof.} It is enough to prove that $G$ contains all
transformations of the following form $$x\mapsto x,\; y\mapsto y,\; z\mapsto z+aM,\;
t\mapsto t,\;\; a\in K,$$ where $M$ is a monomial.

\medskip
{\bf Step 1.} Let $$M=a\prod_{i=1}^m x^{k_i}y^{l_i}\;\;\; \mbox{or}\;\;\;
M=a\prod_{i=1}^m y^{l_0}x^{k_i}y^{l_i}$$ or $$M=a\prod_{i=1}^m
x^{k_i}y^{l_i}\;\;\; \mbox{or}\;\;\; M=a\prod_{i=1}^m
x^{k_i}y^{l_i}x^{k_{m+1}}.$$
Define the height of $M$, $H(M)$, to be the
number of segments comprised of a specific generator - such as $x^k$ - in the word $M$.
(For instance, $H(a\prod_{i=1}^m x^{k_i}y^{l_i}x^{k_{m+1}})=2m+1$.)
Using induction on $H(M)$, one can reduce to the case when
$M=yx^k$. Let $M=M'x^k$ such that $H(M')<H(M)$. (Case when
$M=M'y^l$ is obviously similar.) Let
$$\phi: x\to x,\; y\to y,\; z\to z+M',\; t\to t.$$
$$\alpha: x\to x,\; y\to y,\; z\to z,\; t\to t+zx^k.$$
Then $$\phi^{-1}\circ\alpha\circ\phi: x\to x,\; y\to y,\; z\to z,\; t\to
t-M+ zx^k.$$

The automorphism $\phi^{-1}\circ\alpha\circ\phi$ is the composition of
automorphisms $$\beta: x\to x,\; y\to y,\; z\to z,\; t\to t-M$$ and
$$\gamma: x\to x,\; y\to y,\; z\to z, \;t\to t+zx^k.$$ Observe that $\beta$ is
conjugate to the automorphism $$\beta': x\to x,\; y\to y,\; z\to z-M,\;
t\to t$$ by a linear automorphism $$x\to x,\; y\to y,\; z\to t, \;t\to
z.$$ Similarly, $\gamma$ is conjugate to the automorphism $$\gamma':
x\to x,\; y\to y,\; z\to z+yx^k,\; t\to t.$$ We have thus reduced to the case
when $M=x^k$ or $M=yx^k$.

\medskip
{\bf Step 2.} Consider automorphisms $$\alpha: x\to x,\; y\to y+x^k,\;
z\to z,\; t\to t$$ and $$\beta: x\to x,\; y\to y,\; z\to z,\; t\to t+azy.$$
Then $$\alpha^{-1}\circ\beta\circ\alpha: x\to x,\; y\to y,\; z\to z,\;
t\to t+azx^k+azy.$$ It is a composition of the automorphism $$\gamma:
x\to x,\; y\to y,\; z\to z,\; t\to t+azx^k$$ which is conjugate to the needed
automorphism $$\gamma': x\to x,\; y\to y,\; z\to z+yx^k,\; t\to t$$ and an
automorphism $$\delta: x\to x,\; y\to y,\; z\to z,\; t\to t+azy,$$ which is
conjugate to the automorphism $$\delta': x\to x,\; y\to y,\; z\to
z+axy,\; t\to t$$ and then to the automorphism $$\delta'': x\to x,\;
y\to y,\; z\to z+xy,\; t\to t$$ (using similarities). We have reduced the
problem to proving the statement $$G\ni\psi_M,\;\; M=x^k$$ for all $k$.

\medskip
{\bf Step 3.} Obtain the automorphism $$x\to x,\; y\to y+x^n,\; z\to
z,\; t\to t.$$ This problem is similar to the commutative case of $K[x_1,\dots,x_n]$
(cf. Section \ref{ScAutTameCoommN2}).

Proposition \ref{PrGTame2} is proved.

Returning to the general case $n\geq 4$, let us formulate the remark made after Lemma \ref{LmMult1} as
follows:

\begin{lemma}     \label{LmMult1Ass}
 Consider the following $T^{n-1}$ action: $$x_1\to \lambda^I
x_1,\; x_j\to \lambda_j x_j,\;\; j>1; \;\;\lambda^{I}=\lambda_2^{i_2}\cdots\lambda_n^{i_n}.$$ Then the set $S$ of
automorphisms commuting with this action is generated by the following
automorphisms: $$x_1\to x_1+H,\; x_i\to x_i;\;\; i>1,$$ where $H$ is any homogenous
polynomial of total degree $i_2+\cdots+i_n$.
\end{lemma}

Proposition \ref{PrGTame2} and Lemma \ref{LmMult1Ass} imply

\begin{corollary}                  \label{CoLmMult1Ass}
Let $\Psi\in\Aut(\TAut_0(K\langle x_1,\dots,x_n\rangle))$
stabilize all elements of torus and linear automorphisms,
$$\phi_P: x_n\to x_n+P(x_1,\dots,x_{n-1}),\; x_i\to x_i,\;
i=1,\dots,n-1.$$ Let $P=\sum_IP_I$, where $P_I$ is the homogenous component
of $P$ of multi-degree $I$. Then

a) $\Psi(\phi_P): x_n\to x_n+P^\Psi(x_1,\dots,x_{n-1}),\; x_i\to
x_i,\; i=1,\dots,n-1$.

b) $P^\Psi=\sum_I P_I^\Psi$; here $P_I^\Psi$ is homogenous  of
multi-degree $I$.

c) If $I$ has positive degree with respect to one or two variables,
then $P_I^\Psi=P_I$.
\end{corollary}

Let $\Psi\in\Aut(\TAut_0(K\langle x_1,\dots,x_n\rangle))$
stabilize all elements of torus and linear automorphisms,
$$\phi: x_n\to x_n+P(x_1,\dots,x_{n-1}),\; x_i\to x_i,\;
i=1,\dots,n-1.$$

Let $\varphi_Q: x_1\to x_1,\;x_2\to x_2,\; x_i\to x_i+Q_i(x_1,x_2),\;
i=3,\dots,n-1,\; x_n\to x_n;\;$ $Q=(Q_3,\dots,Q_{n-1})$. Then
$\Psi(\varphi_Q)=\varphi_Q$ by Proposition \ref{PrGTame2}.

\begin{lemma}        \label{LmQPsi12}
a) $\varphi_Q^{-1}\circ\phi_P\circ\varphi_Q =\phi_{P_Q}$, where
$$P_Q(x_1,\dots,x_{n-1})=P(x_1,x_2,x_3+Q_3(x_1,x_2),\dots,x_{n-1}+Q_{n-1}(x_1,x_2)).$$

b) Let $P_Q=P_Q^{(1)}+P_Q^{(2)}$, $P_Q^{(1)}$ consist of all
terms containing one of the variables $x_3,\dots,x_{n-1}$, and let
$P_Q^{(1)}$ consist of all terms containing just
$x_1$ and $x_{2}$. Then
$$P^\Psi_Q=P_Q^\Psi=P_Q^{(1)\Psi}+P_Q^{(2)\Psi}=P_Q^{(1)\Psi}+P_Q^{(2)}$$.
\end{lemma}

\begin{lemma}        \label{LmQPSub}
If $P_Q^{(2)}=R_Q^{(2)}$ for all $Q$ then $P=R$.
\end{lemma}

{\bf Proof.} It is enough to prove that if $P\ne 0$ then
$P^{(2)}_Q\ne 0$ for appropriate $Q=(Q_3,\dots,Q_{n-1})$. Let
$m=\deg(P), \;Q_i=x_1^{2^{i+1}m}x_2^{2^{i+1}m}$. Let $\hat{P}$
be the highest-degree component of $P$, then
$\hat{P}(x_1,x_2,Q_3,\dots,Q_{n-1})$ is the highest-degree component of
$P^{(2)}_Q$. It is enough to prove that
$$\hat{P}(x_1,x_2,Q_3,\dots,Q_{n-1})\ne 0.$$ Let $x_1\prec x_2\prec
x_2\prec\cdots\prec x_{n-1}$ be the standard lexicographic order. Consider the lexicographically minimal
term $M$ of $\hat{P}$. It is easy to see that the term
$$M|_{Q_i\to x_i},\; \;i=3,\;n-1$$ cannot cancel with any other term
$$N|_{Q_i\to x_i},\;\; i=3,\;n-1$$ of
$\hat{P}(x_1,x_2,Q_3,\dots,Q_{n-1})$. Therefore
$\hat{P}(x_1,x_2,Q_3,\dots,Q_{n-1})\ne 0$.

Lemmas \ref{LmQPsi12} and \ref{LmQPSub} imply

\begin{corollary}                  \label{CoSemifinal}
Let $\Psi\in\Aut(\TAut_0(K\langle x_1,\dots,x_n\rangle))$
stabilize all elements of torus and linear automorphisms. Then
$P^\Psi=P$, and $\Psi$ stabilizes all elementary automorphisms and therefore the entire group $\TAut_0(K\langle x_1,\dots,x_n\rangle)$.
\end{corollary}

We obtain the following

\begin{proposition}            \label{PrFinalNGEQ4}
Let $n\ge 4$ and let $\Psi\in\Aut(\TAut_0(K\langle
x_1,\dots,x_n\rangle))$ stabilize all elements of torus and
linear automorphisms. 
Then either $\Psi=\Id$ or $\Psi$ acts as conjugation by the mirror anti-automorphism.
\end{proposition}

Let $n\ge 4$. Let $\Psi\in\Aut(\TAut_0(K\langle
x_1,\dots,x_n\rangle))$ stabilize all elements of torus and
linear automorphisms. Denote by $EL$ an elementary automorphism
$$
EL:x_1\to x_1,\;\ldots,\;x_{n-1}\to x_{n-1},\;x_n\to x_n+x_1x_2
$$
(all other elementary automorphisms of this form, i.e. $x_k\to x_k+x_ix_j,\;x_l\to x_l$ for $l\neq k$ and $k\neq i,\;k\neq j,\;i\neq j$, are conjugate to one another by permutations of generators).

We have to prove that $\Psi(EL)=EL$ or
$\Psi(EL): x_i\to x_i;\; i=1,\dots,x_{n-1},\; x_n\to x_n+x_2x_1$. The latter corresponds to $\Psi$ being the conjugation with the mirror anti-automorphism
of $K\langle x_1,\dots,x_n\rangle$.

Define for some $a,b\in K$ $$x*_{a,b}y=axy+byx.$$

Then, in any of the above two cases, $$\Psi(EL): x_i\to
x_i;\; i=1,\dots,x_{n-1},\; x_n\to x_n+x_1*_{a,b}x_2$$
for some $a,b$.

The following lemma is elementary:

\begin{lemma}                     \label{LmAssoc}
The operation $*=*_{a,b}$ is associative if and only if  $ab=0$. 
\end{lemma}

The associator of $x,\;y,\;\text{and}\;z$ is given by
\begin{eqnarray*}
\lbrace x,y,z\rbrace_{*}\equiv (x*y)*z-x*(y*z)=\\
ab(zx-xz)y+aby(xz-zx)=ab[y,[x,z]].
\end{eqnarray*}


Now we are ready to prove Proposition \ref{PrFinalNGEQ4}. For simplicity we treat only the case $n=4$ -- the general case is dealt with analogously. Consider the
automorphisms $$\alpha: x\to x,\; y\to y,\; z\to z+xy,\; t\to t,$$
$$\beta: x\to x,\; y\to y,\; z\to z,\; t\to t+xz,$$ $$h:x\to x,\; y\to y,\;
z\to z,\; t\to t-xz.$$ (Manifestly $h=\beta^{-1}$.)
Then $$\gamma=h\alpha^{-1}\beta\alpha=[\beta,\alpha]: x\to
x,\; y\to y,\; z\to z,\; t\to t-x^2y.$$
Note that $\alpha$ is conjugate to $\beta$ via a generator permutation
$$
\kappa:x\to x,\;y\to z,\;z\to t,\;t\to y,\;\;\kappa\circ\alpha\circ\kappa^{-1}=\beta
$$
and $$\Psi(\gamma): x\to x,\; y\to
y,\; z\to z,\; t\to t-x*(x*y).$$

Let $$\delta: x\to x,\; y\to y,\; z\to z+x^2,\; t\to t,$$
$$\epsilon: x\to
x,\; y\to y,\; z\to z,\; t\to t+zy.$$ Let
$\gamma'=\epsilon^{-1}\delta^{-1}\epsilon\delta$. Then $$\gamma':
x\to x,\; y\to y,\; z\to z,\; t\to t-x^2y.$$ On the other hand we have
$$\varepsilon=\Psi(\epsilon^{-1}\delta^{-1}\epsilon\delta): x\to x,\;
y\to y, \;z\to z,\; t\to t-(x^2)*y.$$ We also have
$\gamma=\gamma'$. Equality $\Psi(\gamma)=\Psi(\gamma')$ is
equivalent to the equality $x*(x*y)=x^2*y$. This implies $x*y=xy$ and we are done.

\subsection{The group $\Aut_{\Ind}(\TAut(K\langle x,y,z\rangle))$}  \label{SbSc3VrbFreeAss}
This
is the most technically loaded part of the present study. At the moment we are unable to
accomplish the objective of describing the entire group
$\Aut\TAut(K\langle x,y,z\rangle)$. In this
section we will determine only its subgroup
 $\Aut_{\Ind}\TAut_0(K\langle
x,y,z\rangle)$, i.e. the group of $\Ind$-scheme automorphisms, and prove
Theorem \ref{ThTAss3Ind}. We use the approximation results of
Section \ref{SbScAprox}. In what follows we suppose  that
$\Ch(K)\ne 2$. As in the preceding chapter, $\{x,y,z\}_*$ denotes the associator of $x,y,z$ with respect
to a fixed binary linear operation $*$, i.e.
$$\{x,y,z\}_*=(x*y)*z-x*(y*z).$$ 

\begin{proposition}                      \label{LeTAss3Ind}
Let $\Psi\in\Aut_{\Ind}(\TAut_0(K\langle x,y,z\rangle))$
stabilize all linear automorphisms. Let $$\phi: x\to x,\; y\to y,\; z\to
z+xy.$$ Then either $$\Psi(\phi): x\to x,\; y\to y,\; z\to z+axy$$ or
$$\Psi(\phi): x\to x,\; y\to y,\; z\to z+byx$$
for some $a,b\in K$.
\end{proposition}

{\bf Proof.}\ Consider the automorphism $$\phi: x\to x,\; y\to y,\; z\to
z+xy.$$ Then $$\Psi(\phi): x\to x,\; y\to y,\; z\to z+x*y,$$ where $x*y=axy+byx$. Let $a\neq 0$.
We can make the star product $*=*_{a,b}$ into $x*y=xy+\lambda yx$ by conjugation with the mirror anti-automorphism and appropriate linear substitution. We therefore need to
prove that $\lambda=0$, which implies $\Psi(\phi)=\phi$.

The following two lemmas are proved by straightforward computation.

\begin{lemma}     \label{LmAssociator}
Let $A=K\langle x,y,z\rangle$. Let $f*g=fg+\lambda fg$. Then
$\{f,g,h\}_*=\lambda[g,[f,h]]$.

In particular $\{f,g,f\}_*=0$,
$f*(f*g)-(f*f)*g=-\{f,f,g\}_*=\lambda [f,[f,g]]$,\\
$(g*f)*f-g*(f*f)=\{g,f,f\}_*=\lambda [f,[f,g]]$.
\end{lemma}


\begin{lemma}     \label{LmxyyzCommutant}
Let $\varphi_1: x\to x+yz,\; y\to y,\; z\to z$; $\varphi_2: x\to x,\;
y\to y,\; z\to z+yx$;
$\varphi=\varphi_2^{-1}\varphi_1^{-1}\varphi_2\varphi_1$. Then
modulo terms of order $\ge 4$ we have:

$$\varphi: x\to x+y^2x,\; y\to y,\; z\to z-y^2z$$ and $$\Psi(\varphi):
x\to x+y*(y*x),\; y\to y,\; z\to z-y*(y*z).$$
\end{lemma}


\begin{lemma}     \label{LmElemStepxxy}
a) Let  $\phi_l: x\to x,\; y\to y,\; z\to z+y^2x$.  Then $$\Psi(\phi_l):
x\to x,\; y\to y,\; z\to z+y*(y*x).$$

b) Let  $\phi_r: x\to x,\; y\to y,\; z\to z+xy^2$.  Then $$\Psi(\phi_r):
x\to x,\; y\to y,\; z\to z+(x*y)*y.$$

\end{lemma}

{\bf Proof.} According to the results of the previous section we have
$$\Psi(\phi_l): x\to x,\; y\to y,\; z\to z+P(y,x)$$ where $P(y,x)$ is
homogenous of degree 2 with respect to $y$ and degree 1 with
respect to $x$. We have to prove that $H(y,x)=P(y,x)-y*(y*x)=0$.

Let $\tau: x\to z,\; y\to y,\; z\to x;\; \tau=\tau^{-1},\;\;
\phi'=\tau\phi_l\tau^{-1}: x\to x+y^2z,\; y\to y,\; z\to z$. Then
$\Psi(\phi_l'): x\to x+P(y,z),\; y\to y,\; z\to z$.

Let $\phi_l''=\phi_l\phi_l': x\to x+P(y,z),\; y\to y,\; z\to z+P(y,x)$
modulo terms of degree $\ge 4$.

Let $\tau: x\to x-z,\; y\to y,\; z\to z$ and let $\varphi_2$, $\varphi$~ be the
automorphisms described in Lemma \ref{LmxyyzCommutant}.

Then $$T=\tau^{-1}\phi_l^{-1}\tau\phi_l'':x\to x,\;y\to y,\;z\to z$$ modulo
terms of order $\ge 4$.

On the other hand $$\Psi(T): x\to x+H(y,z)-H(y,x),\; y\to y,\; z\to
z+P$$ modulo terms of order $\ge 4$. Because
$\deg_y(H(y,x)=2,\;\deg_x(H(y,x))=1$ we get $H=0$.

 Proof of  b) is  similar.

\begin{lemma}     \label{LmElemSquare}
a)  Let  $$\psi_1: x\to x+y^2,\; y\to y,\; z\to z;\;\; \psi_2: x\to x,\;
y\to y,\; z\to z+x^2.$$  Then
$$[\psi_1,\psi_2]=\psi_2^{-1}\psi_1^{-1}\psi_2\psi_1: x\to x,\; y\to
y,\; z\to z+y^2x+xy^2,$$ $$\Psi([\psi_1,\psi_2]): x\to x,\;y\to y,\;
z\to z+(y*y)*x+x*(y*y).$$

b) $$\phi_l^{-1}\phi_r^{-1}[\psi_1,\psi_2]: x\to x,\;y\to y,\;z\to z$$
modulo terms of order $\ge 4$ but
\begin{gather*}
\Psi\left(\phi_l^{-1}\phi_r^{-1}[\psi_1,\psi_2]\right): x\to x,\;y\to y,\\
z\to
z+(y*y)*x+x*(y*y)-(x*y)*y-y*(y*x)=\\
=z+4\lambda [x[x,y]]
\end{gather*}
modulo terms of order $\ge 4$.
\end{lemma}

{\bf Proof.} a) can be obtained by direct computation. b) follows
from a) and the lemma \ref{LmAssociator}.


Proposition \ref{LeTAss3Ind} follows from Lemma
\ref{LmElemSquare}.

We need a few auxiliary lemmas. The first one is an analogue of
the hiking procedure from
\cite{BelovUzyRowenSerdicaBachtur,BelovIAN}.

\begin{lemma}                  \label{LeHiking}
Let $K$ be  algebraically closed, and let \ $n_1,\dots,n_m$
be positive integers. Then there exist  $k_1,\dots,k_s\in {\mathbb Z}$ and
$\lambda_1,\dots,\lambda_s\in K$ such that

\begin{itemize}
    \item $\sum k_i=1$ modulo $\Ch(K)$ (if $\Ch(K)=0$ then $\sum
    k_i=1$).
    \item $\sum_i k_i^{n_j}\lambda_i=0$ for all $j=1,\dots,m$.
\end{itemize}
\end{lemma}

For $\lambda\in K$ we define an automorphism $\psi_\lambda: x\to
x,\; y\to y,\; z\to\lambda z$.

The next lemma provides for some translation between the language of
polynomials and the group action language. It is similar to the hiking
process \cite{BelovIAN,BelovUzyRowenSerdicaBachtur}.

\begin{lemma}         \label{LeForHiking}
Let $\varphi\in K\langle x,y,z\rangle$. Let $\varphi(x)=x,\;
\varphi(y)=y+\sum_i R_i+R',\; \varphi(z)=z+Q$. Let $\deg(R_i)=N$, let also the
degree of all monomials in $R'$ be greater than $N$, and let the degree
of all monomials in $Q$ be greater than or equal to $N$. Finally, assume $\deg_z(R_i)=i$ and the
$z$-degree of all monomials of $R_1$ greater than $0$.

Then

a) $\psi_\lambda^{-1}\varphi\psi_\lambda: x\to x,\; y\to y+\sum_i
\lambda^iR_i+R'',\; z\to z+Q'$. Also the total degree of all monomials comprising $R'$ is
greater than $N$, and the degree of all monomials of $Q$ is greater than or
equal to $N$.

b) Let $\phi=\prod
\left(\psi_{\lambda^{-1}_i}\varphi\psi_{\lambda_i}\right)^{k_i}$.
Then

$$\phi: x\to x,\; y\to y+\sum_i R_i\lambda_i^{k_i}+S,\; z\to z+T$$
where the degree of all monomials of $S$ is greater than $N$ and the degree of
all monomials of $T$ is greater than or equal to $N$.
\end{lemma}

{\bf Proof.} a) By direct computation. b) is a consequence of a).

{\bf Remark.} In the case of characteristic zero, the condition of $K$
being algebraically closed can be dropped. After hiking for several steps, we need to prove just

\begin{lemma}                  \label{LeHikingCh0}
Let $\Ch(K)=0$, let $n$ be a positive integer. Then there exist
$k_1,\dots,k_s\in {\mathbb Z}$ and
$\lambda_1,\dots,\lambda_s\in K$ such that

\begin{itemize}
    \item $\sum k_i=1$.
    \item $\sum_i k_i^{n}\lambda_i=0$.
\end{itemize}
\end{lemma}

Using this lemma we can cancel out all the terms in the product in the
Lemma \ref{LeForHiking} except for the constant one. The proof of Lemma
\ref{LeHikingCh0} for any field of zero characteristic can be
obtained through the following observation:

\begin{lemma}     \label{LeInclExcl}
$$\left(\sum_{i=1}^n \lambda_i\right)^n-\sum_j
\left(\lambda_1+\cdots+\widehat{\lambda_j}+\cdots+\lambda_n\right)^n+\cdots+
$$
$$+(-1)^{n-k}\sum_{i_1<\cdots<i_k}\left(x_{i_1}+\cdots+x_{i_k}\right)^n+
\cdots+(-1)^{n-1}\left(x_1^n+\cdots+x_n^n\right)=n!\prod_{i=1}^nx_i
$$
and if $m<n$ then
$$\left(\sum_{i=1}^n \lambda_i\right)^m-\sum_j
\left(\lambda_1+\cdots+\widehat{\lambda_j}+\cdots+\lambda_n\right)^m+\cdots+
$$
$$+(-1)^{n-k}\sum_{i_1<\cdots<i_k}\left(x_{i_1}+\cdots+x_{i_k}\right)^m+
\cdots+(-1)^{n-1}\left(x_1^m+\cdots+x_n^m\right)=0.
$$
\end{lemma}

The lemma \ref{LeInclExcl} allows us to replace the $n$-th powers by
product of constants, after that the statement of Lemma
\ref{LeHikingCh0} becomes transparent.

\begin{lemma}               \label{Le2xyAprox}
Let $\varphi: x\to x+R_1,\; y\to y+R_2,\; z\to z'$, such that the total degree of all
monomials in $R_1,\; R_2$ is greater than or equal to $N$. Then for $\Psi(\varphi):
x\to x+R_1',\; y\to y+R_2',\; z\to z''$ with the total degree of all monomials in
$R_1', R_2'$ also greater than or equal to $N$.
\end{lemma}

{\bf Proof.} Similar to the proof of Theorem
\ref{ThMainTechn}.

Lemmas \ref{Le2xyAprox}, \ref{LeForHiking}, \ref{LeHiking} imply
the  following statement.

\begin{lemma}                  \label{LmTechInvar}
Let $\varphi_j\in \Aut_0(K\langle x,y,z\rangle),\; j=1,2$, such that
$$\varphi_j(x)=x,\ \varphi_j(y)=y+\sum_i R_i^j+R'_j,\;
\varphi_j(z)=z+Q_j.$$ Let $\deg(R_i^j)=N$, and suppose that the degree of all monomials
in $R'_j$ is greater than $N$, while the degree of all monomials in $Q$ is
greater than or equal to $N$; $\deg_z(R_i)=i$, and the $z$-degree of all monomials in
$R_1$ is positive. Let $R_0^1=0,\; R_0^2\ne 0$.

Then $\Psi(\varphi_1)\ne\varphi_2$.
\end{lemma}

Consider the automorphism

$$\phi: x\to x,\; y\to y,\; z\to z+P(x,y).$$
Let $\Psi\in\Aut_{\Ind}\TAut_0(k\langle x,y,z\rangle)$ stabilize the
standard torus action pointwise. Then $$\Psi(\phi): x\to x,\; y\to y,\; z\to
z+Q(x,y).$$ We denote
$$\bar{\Psi}(P)=Q.$$ Our goal is to prove that $\bar{\Psi}(P)=P$
for all $P$ if $\Psi$ stabilizes all linear automorphisms and
$\bar{\Psi}(xy)=xy$. We proceed by strong induction on total degree. The base case corresponds to $k=1$ and $l=1$ and is assumed. We then heave

\begin{lemma}
$$\bar{\Psi}(x^ky^l)=x^ky^l$$
provided that $\bar{\Psi}(P)=P$ for all monomials $P(x,y)$ of total degree $< k+l$.
\end{lemma}

{\bf Proof.}

Let $$\phi: x\to x,\; y\to y,\; z\to z+x^ky^l,$$
$$\varphi_1: x\to x+y^l,\; y\to y,\; z\to z,$$ $$\varphi_2: x\to x,\; y\to
y+x^k,\; z\to z,$$ $$\varphi_3: x\to x,\; y\to y,\; z\to z+xy,$$ $$h:
x\to x,\; y\to y,\; z\to z-x^{k+1}.$$ Then, for $k>1$ and $l>1$

\begin{gather*}
g=h\varphi_3^{-1}\varphi_1^{-1}\varphi_2^{-1}\varphi_3\varphi_1\varphi_2:\\
x\mapsto x-y^l+(y-(x-y^l)^k)^l,\\
y\mapsto y-(x-y^l)^k+(x-y^l+(y-(x-y^l)^k)^l)^k,\\
z\mapsto z-xy-x^{k+1}+(x-y^l)(y-(x-y^l)^k).
\end{gather*}
Observe that the height of $g(x)-x$, $g(y)-y$ and $g(z)-z$ is at least $k+l-1$, when $k>1$ or $l>1$.
We then use Theorem
\ref{ThMainTechn} and the induction step. Applying $\Psi$ yields the result because
$\Psi(\varphi_i)=\varphi_i,\; i=1,2,3$ and $\varphi(H_N)\subseteq
H_N$ for all $N$. The lemma is proved.

Let $$M_{k_1,\dots,k_s}=x^{k_1}y^{k_2}\cdots y^{k_s}$$ for even
$s$ and
$$M_{k_1,\dots,k_s}=x^{k_1}y^{k_2}\cdots x^{k_s}$$ for odd
$s$, $k=\sum_{i=1}^n k_i$. Then
$$M_{k_1,\dots,k_s}=M_{k_1,\dots,k_{s-1}}y^{k_s}$$ for even $s$
and
$$M_{k_1,\dots,k_s}=M_{k_1,\dots,k_{s-1}}x^{k_s}$$ for odd
$s$.

We have to prove that
$\bar{\Psi}(M_{k_1,\dots,k_{s}})=M_{k_1,\dots,k_{s}}$.
By induction we may assume that
$\bar{\Psi}(M_{k_1,\dots,k_{s-1}})=M_{k_1,\dots,k_{s-1}}$.

For any monomial $M=M(x,y)$ we define an automorphism
$$\varphi_M: x\to x,\; y\to y,\; z\to z+M.$$

We also define the automorphisms $$\phi_k^e: x\to x,\; y\to y+zx^k,\; z\to
z$$ and $$\phi_k^o: x\to x+zy^k,\; y\to y,\; z\to z.$$ We will present the
case of even $s$ - the odd $s$ case is similar.

Let $D_{zx^k}^e$ be a derivation of $K\langle x,y,z\rangle$ such
that $D_{zx^k}^e(x)=0$,
 $D_{zx^k}^e(y)=zx^k$, $D_{zx^k}^e(z)=0$.
Similarly, let $D_{zy^k}^o$ be a derivation of $K\langle x,y,z\rangle$
such that $D_{zy^k}^o(y)=0$, $D_{zx^k}^o(x)=zy^k$,
$D_{zy^k}(z)^o=0$.

The following lemma is proved by direct computation:

\begin{lemma}
Let
$$u=\phi_{k_{s}}^e{}^{-1}\varphi(M_{k_1,\dots,k_{s-1}})^{-1}\phi_{k_{s}}^e\varphi(M_{k_1,\dots,k_{s-1}})$$
for even $s$ and
$$u=\phi_{k_{s}}^o{}^{-1}\varphi(M_{k_1,\dots,k_{s-1}})^{-1}\phi_{k_{s}}^o\varphi(M_{k_1,\dots,k_{s-1}})$$
for odd $s$. Then

$$u:x\to x,\; y\to y+M_{k_1,\dots,k_s}+N',\; z\to
z+D^e_{zx^k}(M_{k_1,\dots,k_{s-1}})+N$$ for even $s$ and
$$u:x\to x+M_{k_1,\dots,k_s}+N',\; y\to y,\; z\to
z+D^o_{zx^k}(M_{k_1,\dots,k_{s-1}})+N$$ for odd $s$,
where $N$, $N'$ are sums of terms of degree
$>k=\sum_{i=1}^sk_i$.
\end{lemma}

Let $\psi(M_{k_1,\dots,k_{s}}): x\to x,\; y\to y,\; z\to
z+M_{k_1,\dots,k_{s}}$, $$\alpha_e: x\to x,\; y\to y-z,\;z\to z,
\alpha_o: x\to x-z,\; y\to y,\;z\to z,$$ Let $P_M=\Psi(M)-M$. Our goal
is to prove that $P_M=0$.

Let
$$v=\psi(M_{k_1,\dots,k_{s}})^{-1}\alpha_{e}\psi(M_{k_1,\dots,k_{s}})u\alpha_{e}^{-1}$$
for even $s$ and
$$v=\psi(M_{k_1,\dots,k_{s}})^{-1}\alpha_{o}\psi(M_{k_1,\dots,k_{s}})u\alpha_{o}^{-1}$$
for odd $s$.

The next  lemma  is also proved by direct computation:

\begin{lemma}      \label{LmFinalType}
a) $$v: x\to x,\; y\to y+H,\; z\to z+H_1+H_2$$ for even $s$ and
$$v: x\to x+H,\; y\to y,\; z\to z+H_1+H_2$$
for odd $s$

b)
$$\Psi(v): x\to x, y\to y+P_{M_{k_1,\dots,k_{s}}}+\widetilde{H}, z\to z+\widetilde{H_1}+\widetilde{H_2}$$
for even $s$ and
$$\Psi(v): x\to x+P_{M_{k_1,\dots,k_{s}}}+\widetilde{H}, y\to y, z\to z+\widetilde{H_1}+\widetilde{H_2}$$
for odd $s$,
where $H_2$, $\widetilde{H_2}$ are sums of terms of degree
greater than $k=\sum_{i=1}^s k_i$, $H$, $\widetilde{H}$
are sums of terms of degree $\ge k$ and positive $z$-degree, $H_1$, $\widetilde{H_1}$ are sums of terms of
degree $k$ and positive $z$-degree.
\end{lemma}

%
%
%
%

{\bf Proof of Theorem \ref{ThTAss3Ind}.} Part b) follows from
part a). In order to prove a) we are going to show that
$\bar{\Psi}(M)=M$ for any monomial $M(x,y)$ and for any
$\Psi\in\Aut_{\Ind}(\TAut(\langle x,y,z\rangle))$ stabilizing the
standard torus action $T^3$ and $\phi$. The automorphism
$\Psi(\Phi_M)$ has the form described in Lemma
\ref{LmFinalType}. But in this case Lemma \ref{LmTechInvar} implies
$\bar{\Psi}(M)-M=0$.

%
%
%
%
%

\section{Some open questions concerning the tame automorphism group}

As the conclusion of the paper, we would like to raise the following questions.

\begin{enumerate}
    \item Is it true that any  automorphism $\varphi$ of $\Aut(K\langle
x_1,\dots,x_n\rangle)$ (in the group-theoretic sense - that is, not necessarily an automorphism preserving the $\Ind$-scheme structure) for
$n=3$ is semi-inner, i.e. is a conjugation by some automorphism
or mirror anti-automorphism?
    \item Is it true that $\Aut(K\langle
x_1,\dots,x_n\rangle)$ is generated by affine automorphisms and
automorphism $x_n\to x_n+x_1x_2,\; x_i\to x_i,\; i\ne n$? For $n\ge 5$
it seems to be easier and the answer is probably positive, however for $n=3$ the
answer is known to be negative, cf. Umirbaev \cite{U} and Drensky and Yu
\cite{DYuStrongAnik}. For $n\ge 4$ we believe the answer is
positive.
    \item  Is it true that $\Aut(K[x_1,\dots,x_n])$ is generated by linear automorphisms and
automorphism $x_n\to x_n+x_1x_2,\; x_i\to x_i,\; i\ne n$? For $n=3$
the answer is negative: see the proof of the Nagata conjecture
\cite{SU1,SU2,UY}. For $n\ge 4$ it is plausible that the answer is
positive.
    \item Is any  automorphism $\varphi$ of $\Aut(K\langle
x,y,z\rangle)$ (in the group-theoretic sense)
semi-inner?
    \item Is it true that the conjugation in Theorems \ref{ThAutTAut}
    and \ref{ThAutTAutFreass} can be done by some tame automorphism?
    Suppose $\psi^{-1}\varphi\psi$ is tame for any
    tame $\varphi$. Does it follow that $\psi$ is tame?
    \item Prove Theorem \ref{ThTAss3Ind} for
    $\Ch(K)=2$. Does it hold on the set-theoretic level, i.e. $\Aut(\TAut(K\langle
x,y,z\rangle))$ are generated by conjugations by an automorphism or the
mirror anti-automorphism?
\end{enumerate}

Similar questions can be formulated for nice automorphisms.

\section{Acknowledgements}

The authors would like to thank to J. P. Furter, T. Kambayashi, R. Lipyanski, B. Plotkin, E. Plotkin, A. Regeta,  and M. Zaidenberg for stimulating discussion.
We are grateful to Eliyahu Rips for having kindly agreed to include and
use his crucial results.

The authors also thank Shanghai University and Shenzhen university
for warm hospitality and stimulating atmosphere during their visits,
when parts of this project were carried out.


\end{document}